\input amstex
\documentstyle{amsppt}
\nologo
\TagsOnRight
\magnification 1200
\voffset -1cm

\def\on#1{\expandafter\def\csname#1\endcsname{{\operatorname{#1}}}} \on{im}
\on{lf} \on{id} \on{dist} \on{deg} \on{Hom}  \on{holink}
\on{eq} \on{Tel} \on{Cyl} \on{Fr}  \on{Int}
\def\fn#1{\expandafter\def\csname#1\endcsname{\operatorname{#1}\def\shortcut
{\if(\next{}\else\,\fi}\futurelet\next\shortcut}} \fn{coker} \fn{mod} \fn{rel}
\let\emb\hookrightarrow \let\imm\looparrowright  \def\Q{\Bbb Q}
\def\R{\Bbb R} \def\?{@!@!@!@!@!} \def\RP{\R\? P} \def\z#1#2{#1\?#2\?}
\def\Z{{\Bbb Z}\def\shortcut{\if/\next{}\z\fi}\futurelet\next\shortcut}
\let\tl\tilde \let\but\setminus \let\x\times \let\eps\varepsilon \on{coind}
\let\tta\vartheta \let\phi\varphi \def\T{_{\?\text{\fiverm T}}} \def\N{\Bbb N}
\def\invlim{\lim\limits_{\longleftarrow}\,} \let\del\vartriangle
\def\dirlim{\lim\limits_{\longrightarrow}\,} \let\imp\Rightarrow
\def\derlim{\lim\limits_{\longleftarrow}{^{\!1}}\,} \on{rk} \fn{lk} \fn{st}
\let\inv\leftrightarrow  \on{mesh} \on{fr}
\def\Cl#1{\overline{#1}}  \def\H{\Cal H}

\topmatter 
\thanks Partially supported by Russian Foundation for Basic Research grants
02-01-00014, 05-01-00993, Mathematics Branch of the Russian Academy of
Sciences Program ``Theoretical Problems of Contemporary Mathematics'' and
President of Russia Grant MK-1844.2005.1.
\endthanks
\address Steklov Mathematical Institute, Division of Geometry and Topology;
Gubkina ~8, Moscow 119991, Russia \endaddress
\email melikhov\@mi.ras.ru; scepin\@mi.ras.ru \endemail

\title The telescope approach to embeddability of compacta \endtitle

\author Sergey A. Melikhov and Evgenij V. Shchepin \endauthor

\abstract
We show that an $n$-dimensional compactum $X$ embeds in $\R^m$, where
$m>\frac{3(n+1)}2$, if and only if $X\x X\but\Delta$ admits an
equivariant map to $S^{m-1}$.
In particular, $X$ embeds in $\R^{2n}$, $n>3$, if and only if the top
power of the (twisted) Euler class of the factor-exchanging involution
on $X\x X\but\Delta$ is trivial.
Assuming that $X$ quasi-embeds in $\R^{2n}$ (i.e\. is an inverse limit of
$n$-polyhedra, embeddable in $\R^{2n}$), this is equivalent to
the vanishing of an obstruction in $\lim^1 H^{2n-1}(K_i)$ over compact
subsets $K_i\i X\x X\but\Delta$.
One application is that an $n$-dimensional ANR embeds in $\R^{2n}$
if it quasi-embeds in $\R^{2n-1}$, $n>3$.

We construct an ANR of dimension $n>1$, quasi-embeddable but not embeddable
in $\R^{2n}$, and an AR of dimension $n>1$, which does not ``movably'' embed
in $\R^{2n}$.
These examples come close to, but don't quite resolve, Borsuk's problem: does
every $n$-dimensional AR embed in $\R^{2n}$?
In the affirmative direction, we show that an $n$-dimensional compactum $X$
embeds in $\R^{2n}$, $n>3$, if $H^n(X)=0$ and $H^{n+1}(X,X\but x)=0$ for
every $x\in X$.

There are applications in the entire metastable range as well.
An $n$-dimensional compactum $X$ with $H^{n-i}(X\but x)=0$ for each $x\in X$
and all $i\le k$ embeds in $\R^{2n-k}$.
This generalizes Bryant and Mio's result that $k$-connected $n$-dimensional
generalized manifolds embed in $\R^{2n-k}$.
Also, an acyclic compactum $X$ embeds in $\R^m$ iff $X\x I$ embeds in
$\R^{m+1}$ iff $X\x$\,(triod) embeds in $\R^{m+2}$.
As a byproduct, we answer a question of T. Banakh on stable embeddability of
the Menger cube.
\endabstract
\endtopmatter

\document

\head 1. Introduction \endhead

It is well-known that all (compact) contractible $n$-polyhedra embed in
$\R^{2n}$ but not all in $\R^{2n-1}$ \cite{We}.
(The case $n=2$ is more subtle: every contractible $2$-polyhedron PL embeds
in a homotopy $4$-sphere, hence by Freedman's work embeds in $\R^4$.)
This paper was motivated by

\proclaim{Problem 1.1 (Borsuk, 1959)} {\rm (see \cite{Bo; IX.2.4})} Does
every $n$-dimensional AR (=contractible, locally contractible compactum)
embed in $\R^{2n}$?
\endproclaim

As noted in \cite{Bo}, the answer is affirmative for $n=1$.
This follows e.g\. from Claytor's theorem (see Theorem 2.1 below).
Contractible $n$-dimensional compacta, non-embeddable in $\R^{2n}$,
were constructed by S. D. Iliadis for $n=1$ \cite{SS} and in \cite{RSS}
for all $n$.
The construction of \cite{RSS} will be discussed in Example 2.17, and here
is another one.

\example{Example 1.2 (the $p$-adic tree)}
Let $X=\Tel(\Z_p@>\lfloor/p\rfloor>>\Z_p@>\lfloor/p\rfloor>>\dots)^*$,
the one-point compactification of the mapping telescope of $p$-fold
self-coverings of the Cantor set, realized here as integral divisions by $p$
in the topological ring $\Z_p$ of $p$-adic integers.
Clearly, $X$ is contractible.
It does not embed in $\R^2$ since it contains a copy of the compactum
$T\x\N^*$, where $T$ is the triod $pt*(\Z/3)$, and $\N^*$ denotes
the one-point compactification of an infinite sequence.
It is shown in Example 2.17 that $T\x\N^*$ does not embed in $\R^2$.
Note that $X$ is the inverse limit of the finite mapping telescopes
$X_n:=\Tel(\Z/p^n@>\lfloor/p\rfloor>>\dots@>\lfloor/p\rfloor>>\Z/p@>>>0)$
and maps $X_n\to X_{n-1}$ shrinking the last mapping cylinder
$\Cyl(\Z/p\to 0)$ to the point $0$ and restricting to the trivial
$p$-covering on its complement.
\endexample

Our initial result was that $n$-dimensional contractible compacta
{\it quasi-embed} in $\R^{2n}$ for $n>2$, that is, are inverse limits of
$n$-polyhedra, embeddable in $\R^{2n}$.
(An equivalent definition is that for each $\eps>0$, there exists
an $\eps$-map $X\to\R^{2n}$, i.e\. a map whose every point-inverse has
diameter $<\eps$.)

\proclaim{Theorem 1.3}
Every $n$-dimensional compactum $X$, $n>2$, such that the map
$H^n(X,X\but x)\to H^n(X)$ is onto for some $x\in X$, quasi-embeds in
$\R^{2n}$.
\endproclaim

This follows from Theorem 2.4 and Lemma 3.6(a).
Note that not every $n$-dimensional contractible compactum can be decomposed
into an inverse limit of contractible (or even acyclic) $n$-polyhedra
--- even if it is an inverse limit of PL balls of some dimension \cite{Ka}.

For ANRs (=locally contractible compacta), Theorem 1.3 holds if
the coefficients are reduced $\bmod 2$; yet there exists an $n$-dimensional
compactum $X$ with $H^n(X;\,\Z/2)=0$, non-quasi-embeddable in $\R^{2n}$ (see
Example 2.6 and the subsequent remark).
This shows, in particular, that not every $n$-dimensional compactum $X$ with
$H^n(X;\,\Z/2)=0$ is an inverse limit of $n$-polyhedra $P_i$ with
$H^n(P_i;\,\Z/2)=0$.

In pursuing further our approach to Problem 1.1, we gradually realized that
we need no less than a complete cohomological obstruction to embeddability of
an arbitrary $n$-dimensional compactum $X$ into $\R^{2n}$ to start with.
This obstruction $\theta(X)$ is the Euler class of the vector bundle
$\tl X\x_{\Z/2}\R^{2n}\to\tl X/(\Z/2)$, where $\tl X=X\x X\but\Delta$
(with the factor exchanging involution), and $\R^{2n}$ is endowed
with the sign action of $\Z/2$ (see other definitions in \S2).
While this is the straightforward generalization of the classical van Kampen
obstruction to embeddability of a compact polyhedron into $\R^{2n}$, from
the geometric viewpoint it involves an additional ``phantom'' term due to
more complex local structure of the compactum $X$, as measured by
non-collarable behavior of $X\x X\but\Delta$ at infinity.

The basic idea behind the completeness of $\theta(X)$ for $n>3$
(Theorem 2.2) is that $X$ embeds into $\R^m$ if and only if an infinite
polyhedron $T$ of a certain kind (namely the mapping telescope of some
inverse sequence of nerves of $X$), which is endowed with a proper
``control'' map to $[0,\infty)$, admits a level-preserving embedding into
$\R^m\x [0,\infty)$.
This was proved in a 1984 paper by the second author and M. A. Shtan'ko
\cite{SS} (not translated into English) using the techniques of embedding
dimension to establish the ``only if'' part in codimension three.
While we only need the easy ``if'' part here (see Criterion 3.2), the
challenge is to construct a telescope $T$ with vanishing cohomological
embedding obstruction, knowing only the vanishing of the cohomological
embedding obstruction for $X$.
To this end we need some algebraic technique to capture the behavior of
$X\x X\but\Delta$ at infinity by means of the extra dimension occurring
in $T$.
Luckily, much of it turned out to be available, in a disguised form, from
the first author's work on isotopic realization \cite{M1}, \cite{M2},
\cite{M3}, where the role of ``local wildness'' of the compactum was played
by ``global wildness'' of the double point set of a continuous map between
polyhedra.
\medskip

Now that the Borsuk problem has been reduced (in dimensions $\ne 2,3$),
to a question of constructing an AR whose deleted product satisfies
a purely cohomological condition (see also Lemma 2.18), it looks even more
formidable.
For this realizability question is not an easy one!

\proclaim{Theorem 1.4} Every $n$-dimensional compactum $X$, $n>3$, such that
the map $H^n(X,X\but x)\to H^n(X)$ is onto for some $x\in X$ and
$H^{n+1}(X,X\but x)=0$ for each $x\in X$, embeds in $\R^{2n}$.
\endproclaim

This follows from Theorem 2.2 and Lemma 3.6.

The reader should not be surprised at seeing the $(n+1)$-dimensional
cohomology of the $n$-dimensional compactum with support in a point.
(Relative cohomology modulo an open set is a counterintuitive object.)
The simplest example of a point where such a group is nontrivial occurs as
the point at infinity in the one-point compactification
$\Tel(P_1\to P_2\to\dots)^*$ of the mapping telescope of a direct sequence
of PL maps between compact $(n-1)$-polyhedra, provided that each
$H^{n-1}(P_{i+1})\to H^{n-1}(P_i)$ has a nonzero cokernel in the
torsion-free part (see Example 2.8).
Coincidentally, this compactum is an AR, as is every
$\Tel(P_1\to P_2\to\dots)^*$ where the $P_i$'s are ANRs.

\proclaim{Problem 1.5} Does every
$\Tel(P_1\to P_2\to\dots)^*$ embed in $\R^{2n}$, where $P_i$ are
$(n-1)$-dimensional polyhedra (or ANRs)?
\endproclaim

A good reason to focus on ARs with just one non-polyhedral point is that
a contractible compactum embeds in $\R^{2n}$ whenever it immerses there,
i.e\. admits a map into $\R^{2n}$ that is injective on a neighborhood of
every point (Corollary 4.9).

At a first glance, direct telescopes of graphs ($n=2$) appear promising.
By a well-known result of Conway--Gordon and Sachs, no matter how
the complete graph $K_6$ is embedded into $\R^3$, some pair of disjoint
cycles in it will be linked with an odd linking number.
(In fact, a given $n$-polyhedron $P$ admits an embedding into $\R^{2n+1}$
where no two disjoint $n$-spheres in the image are linked if and only if
an odd-dimensional analogue $\eta(P)$ of the van Kampen obstruction
vanishes \cite{M5}.)
If $\Gamma\to K_6$ is a $3$-fold covering, say, and its mapping cylinder
has been embedded into $\R^3\x I$ in a level-preserving fashion, then
$\Gamma$ must have a pair of disjoint cycles linked in $\R^3\x\{0\}$ with
linking number equal to $3$ times the linking number in $\R^3\x\{1\}$ of
some pair of disjoint cycles of $K_6$.
Now if we ignore covering theory for a second and assume that there exists
a $3$-fold covering $f\:K_6\to K_6$, then $\Tel(K_6@>f>>K_6@>f>>\dots)$
cannot be embedded into $\R^3\x [0,\infty)$ in a level-preserving fashion,
for the very first $K_6$ would have to possess, for every $n$, a pair of
disjoint cycles linked in $\R^3\x\{0\}$ with linking number divisible by
$3^n$ --- but $K_6$ only contains $10$ pairs of disjoint cycles!
Thus every $3$-fold covering $K_6\to K_6$ would yield a solution to
the Borsuk problem, for the level-preserving condition is not essential here
due to the cohomological nature of the argument.

But, of course, no such coverings exist.
Since a $p$-fold cover of a graph with Euler characteristic $\chi$ has
Euler characteristic $p\chi$, there even exists no direct sequence
$\Gamma_1\to\Gamma_2\to\dots$ of non-trivial coverings of finite graphs with
more than one cycle in at least one connected component.
(This also has implications for direct sequences of branched coverings
between higher-dimensional polyhedra, which one may restrict to maps between
their intrinsic $1$-skeleta, or inductively restrict to branched coverings
between links of the vertices in the intrinsic $0$-skeleta.)

\proclaim{Theorem 1.6 \cite{M6}} Let
$X=\Tel(\Gamma_1@>f_1>>\Gamma_2@>f_2>>\dots)^*$
where $\Gamma_i$ are finite graphs.

(a) If each $f_i$ is onto a finite index subgroup on the $\pi_1$ level,
$X$ embeds into a $Q$ such $Q\x\R^2$ is homeomorphic to $\R^6$.

(b) $X$ embeds into $\R^4$ after amending each $f_i$ by a homotopy.
\endproclaim

The proof is modelled on the well-known proofs that $2$-polyhedra,
$3$-deformable to a point, embed in $\R^4$ and that $2$-dimensional
CW-complexes embed in $\R^4$.

\proclaim{Problem 1.7} Does there exist an $n$-dimensional ANR,
non-embeddable in $\R^{2n}$, which can cover itself with finitely many branch
points?
\endproclaim

Note that the ``intrinsic $1$-skeleton'' of an ANR (whatever it means)
does not need to be an ANR, and so can have an infinitely generated
fundamental group.
\medskip

\subhead Organization of the paper \endsubhead
Results in the double dimension ($n\emb 2n$) are stated and illustrated by
examples in \S2.
(Some further discussion is at the end of \S3.)
Theorems 2.2 and 2.4 are proved in \S3.
Results valid in the entire metastable range are stated and proved in \S4.

\subhead Acknowledgements\endsubhead
The first author is grateful to M. Bestvina and A. N. Dranishnikov for
inviting him to discuss this work with them at their Universities and
their support and useful comments.
We also enjoyed critical and creative comments from P. M. Akhmetiev (in
particular, Example 2.15 is his), stimulating questions of V. M. Buchstaber
and T. Banakh (see \S4), and information provided by A. V. Chernavsky,
J. Keesling and A. B. Skopenkov.
Some results of this paper were presented at Borsuk's 100th Anniversary
Conference (B\c edlewo, July 2005).

\head 2. Results and examples in the double dimension \endhead

Throughout the paper $\bar X$ will denote the quotient of
$\tl X:=X\x X\but\Delta_X$ by $\Z/2=\left<t\mid t^2\right>$, acting on
$\tl X$ by the factor exchanging involution $(x,y)\inv (y,x)$.

The definition of the classical polyhedral van Kampen obstruction $\tta(X)$
(see \cite{M5}) makes perfect sense if $X$ is an $n$-dimensional (metric)
compactum, and then the same argument as in the polyhedral case shows that if
$X$ embeds in $\R^{2n}$ then $\tta(X)=0$.
In more detail, we may define $\tta(X)\in H^{2n}(\bar X)$ to be e.g\.
$\bar g^*(\xi)$, where $g\:X\emb\R^{2n+1}$ is any embedding (which exists by
the classical Menger--N\"obeling--Pontryagin theorem),
$\xi\in H^{2n}(\RP^\infty)\cong\Z/2$ is the generator, and
$\bar g\:\bar X\to\RP^{2n}$ is defined by
$\{x,y\}\mapsto\left<g(x)-g(y)\right>$.
Since any two embeddings $X\emb\R^{2n+2}$ are isotopic%
\footnote{Note that they need not be ambient isotopic, even if $X=I^n$.
On the other hand, it is not hard to see that any two tame embeddings
$X\emb\R^{2n+2}$ are ambient isotopic, whereas Edwards proved that any
embedding in codimension $\ge 3$ can be obtained by an (ambient)
pseudo-isotopy of a tame embedding.}%
, which is proved analogously to the Menger--N\"obeling--Pontryagin theorem,
this well-defines $\tta(X)$.

\proclaim{Theorem 2.1 (Claytor--Skopenkov)} \cite{C}, \cite{Sk1}
A locally connected $1$-dimensi\-onal compactum $X$ embeds in $\R^2$ iff
$\tta(X)=0$.
\endproclaim

\remark{Remark} The idea of proof is as follows.
Claytor proved that every locally connected compactum, non-embeddable in
the plane, contains a copy of either $K_5$ or $K_{3,3}$ (the Kuratowski
graphs) or $C_5$ or $C_{3,3}$ \cite{C}.
Here $C_5$ is the wedge of an arc (with basepoint at an endpoint) and the
one-point compactification of $\hat K_5/T$ (with basepoint at infinity),
where $\hat K_5$ is the infinite cyclic cover of $K_5$ corresponding to
a map $f\:K_5\to S^1$ that contracts all the edges but one to
$0\in S^1=\R/\Z$, and $T$ is an involution on $\hat K_5$ descending to
the reflection $x\mapsto\frac12-x$ on $\R$.
$C_{3,3}$ can be described analogously starting from $K_{3,3}$.
Skopenkov observed that neither $\tl C_5$ nor $\tl C_{3,3}$ admits an
equivariant map to $S^1$ \cite{Sk1}.
(The reader may want to reprove this fact along the lines of Example 2.17
below.)
Finally, by the standard obstruction theory (cf\. \cite{M5} for
the polyhedral case), existence of an equivariant map $\tl X\to S^1$ is
equivalent to the vanishing of $\theta(X)$.
\endremark

\remark{Remark}
The assumption of local connectedness is necessary in Theorem 2.1.
The $p$-adic solenoid $\Sigma_p$, i.e\. the inverse limit of
$\cdots@>p>>S^1@>p>>S^1$ does not embed in the plane.
For, assuming the contrary, we have by the Alexander duality
$H_0(\R^2\but\Sigma_p)\simeq H^1(\Sigma_p)\simeq\Z[\frac1p]$, but the
$0$-homology of an open manifold must be free abelian.
On the other hand, it was noticed by Skopenkov \cite{Sk1} that there exists
an equivariant map $\tl\Sigma_3\to S^1$; therefore $\tta(\Sigma_3)=0$.
\endremark

\proclaim{Theorem 2.2} An $n$-dimensional compactum $X$ embeds in $\R^{2n}$,
$n>3$, if and only if $\tta(X)=0\in H^{2n}(\bar X)$.
\endproclaim

This follows from Criterion 3.2, Lemma 3.3 and Proposition 3.4.

The simplicity of the statement is deceptive.
Let us represent $\bar X$ as a union of compact subsets
$K_0\i K_1\i K_2\i\cdots$.
We have%
\footnote{The {\it derived limit} $\lim^1$ of an inverse sequence of abelian
groups $\{G_i;\pi_i\}$ is the derived functor of the inverse limit $\lim$.
A resolution due to Roos identifies $\lim$ and $\lim^1$ respectively with
the kernel and the cokernel of $\tau\:\prod G_i\to\prod G_i$, sending
$(a_0,a_1,\dots)$ to $(a_0-\pi_1 a_1,a_1-\pi_2 a_2,\dots)$.}
Milnor's exact sequence (see e.g. proofs of Lemmas 2.16 and 2.18 below)
$$0\to\derlim H^{2n-1}(K_i)\to H^{2n}(\bar X)\to\invlim H^{2n}(K_i)\to 0.$$
If the coefficients are reduced $\bmod2$, the derived limit on the left hand
side will trivialize, provided that $X$ is an ANR.
Indeed, by local acyclicity the images of the bonding maps may be assumed
finitely generated \cite{Br}, whereas any inverse sequence of finite groups
satisfies the Mittag-Leffler condition.%
\footnote{An inverse sequence of groups $G_i$ satisfies the
{\it Mittag-Leffler condition}, if for each $i$ there exists a $j>i$
such that for each $k>j$ the image of $G_k$ in $G_i$ equals that of $G_j$.
The ML condition is sufficient for the vanishing of $\lim^1$ and, if all
$G_i$'s are countable, also necessary \cite{Gr}.}
In particular, if $\tta(X)$ lies in this derived limit, it is infinitely
$2$-divisible (note also it always is an element of order two by definition).

\example{Example 2.3 (lim$^{\bold 1}$)}
The simplest example of an inverse sequence of abelian groups with
nontrivial derived limit is the ``$p$-tower'' $\dots@>p>>\Z@>p>>\Z$,
whose derived limit is $\Z_p/\Z$, where $\Z_p$ denotes the $p$-adic integers.
(See \cite{M2; Example 4} for several geometric interpretations of this
computation, involving the $0$-homology of the $p$-adic solenoid.)
Now $\Z_p/\Z$ is the direct sum of an uncountable torsion-free group and
$\Z_{(p)}/\Z$, where $\Z_{(p)}\i\Q$ denotes the localization.
Thus it contains no elements of order two when $p=2$ and precisely one such
element when $p$ is odd.
Another standard example of an inverse sequence with nontrivial $\lim^1$ is
``Jacob's ladder'' $\dots@>\text{incl}>>\bigoplus_{i=1}^\infty\Z
@>\text{incl}>>\bigoplus_{i=0}^\infty\Z$.
The derived limit would vanish (even though the Mittag-Leffler condition
would still fail) if the sums were replaced by products.%
\footnote{In addition to this assertion, the reader who would like to
familiarize herself with $\lim^1$ is invited to directly verify its
non-vanishing for Jacob's ladder and the $p$-tower.
Some references for $\lim^1$ can be found in \cite{M3}.}
\endexample

For brevity we shall use the notation $$\hat H^k(U)=\invlim H^k(K_i),$$
which is supposed to reflect the duality with \v{C}ech homology $\check H_*$.
The image $$\hat\tta(X)\in\hat H^{2n}(\bar X)$$ of $\tta(X)$ obstructs only
quasi-embeddability of $X$ into $\R^{2n}$.
Indeed, let us take each $K_i$ to consist of all unordered pairs
$\{x,y\}$ with $\dist(x,y)\ge 1/i$.
Then a $1/i$-map $f_i\:X\to\R^{2n}$ yields a map
$\bar f_i\:K_i\to\RP^{2n-1}$, $\{x,y\}\mapsto\left<f(x)-f(y)\right>$,
and the projection of $\hat\tta(X)$ to $H^{2n}(K_i)$ equals
$\bar f_i^*(\xi)$, which is zero.

The proof of Theorem 2.2 will be modelled after the (much easier) proof of
Proposition 3.1, which together with completeness of $\tta(X)$ for
compact $n$-polyhedra, $n>2$ (cf\. \cite{M5}) implies

\proclaim{Theorem 2.4} An $n$-dimensional compactum $X$ quasi-embeds in
$\R^{2n}$, $n>2$, if and only if $\hat\tta(X)=0\in\hat H^{2n}(\bar X)$.
\endproclaim

The restriction $n>2$ is known to be necessary already in the polyhedral
case.
In the proof of Theorem 2.2, we encounter the stronger restriction $n>3$
twice, for two apparently distinct reasons.

\proclaim{Problem 2.5} Is $\tta(X)$ complete for $3$-dimensional compacta?
\endproclaim

It follows from the definition that $\tta(X)$ equals $e(\eta_X)^{2n}$, the
top power of the the twisted Euler class of the $2$-cover
$\eta_X\:\tl X\to\bar X$ (see \cite{M5}).
Even though it is an element of order two, a care must be taken to
distinguish it from $w_1(\eta_X)^{2n}$, the top power of the first
Stiefel--Whitney class.
An $n$-polyhedron $P$ with $w_1(\eta_P)^{2n}=0$ but $e(\eta_P)^{2n}\ne 0$
was constructed in \cite{M5} for each $n>1$.
A minor modification of this construction yields

\example{Example 2.6}
There exists an $n$-dimensional compactum $X$ with $H^n(X;\Z/2)=0$,
which does not quasi-embed in $\R^{2n}$.
Specifically, such an $X$ can be obtained from the $n$-skeleton $K$ of
the $(2n+2)$-simplex by replacing one $n$-ball in the interior of each
$n$-simplex $\sigma_i$ with the one-point compactification $T_i^*$ of
the telescope of an inverse sequence of degree two maps
$S^{n-1}\to S^{n-1}$.
(Compare Example 2.12 below, where a {\it direct} telescope of
{\it odd} degree maps is used.)
Since every cocycle with support in a point can be pushed towards
the point at infinity $*_i$ of one of $T_i^*$, where it becomes divisible
by $2$, we have $H^n(X;\Z/2)=0$.

On the other hand, if $g\:X\to\R^{2n}$ is an $\eps$-map with $\eps$
sufficiently small, then $g(T_i^*)\cap g(\sigma_j)=\emptyset$ whenever
$i\ne j$ (where $\sigma_j$ denotes the modified simplex).
Let $Z_i$ be the union of $\sigma_j$'s, disjoint from $\sigma_i$; this
is an $n$-sphere with some small balls replaced by $T_j$'s.
Given a ball $B$ around $*_i$, small enough to be disjoint from
$g(\sigma_j)$'s, the sphere $S^{n-1}_N$ at some level $N$ in the telescope
$T_i$ is contained in $B$.
Then
$\lk(S^{n-1}_0,Z_i)=\frac1{2^N}\lk(S^{n-1}_N,Z_i)=\frac0{2^N}=0\in
H^n(Z_i)\simeq\Z[\frac12]$.
This linking number is the pullback of the generator under the map
$H^{2n-1}(S^{2n-1})\to H^{2n-1}(S^{n-1}_0\x Z_i)$, induced by the restriction
of $\tl g_i$.
It follows (cf\. \cite{M5; proof of Lemma 2.5}) that $\tl g$ (restricted
to $\tl X\but N\Delta$, where $N\Delta$ is a neighborhood of the diagonal
corresponding to $\eps$) is equivariantly homotopic to a map that factors
through $\tl X_i\but N\Delta$, where $X_i$ is constructed similarly to $X$
except that the simplex $\sigma_i$ is left unaltered.
Arguing by induction, we obtain an equivariant map $\tl K\to S^{2n-1}$, which
cannot be.
\endexample

\remark{Remark} If $X$ is a locally acyclic $n$-dimensional compactum such
that the map $H^n(X,X\but x;\,\Z/2)\to H^n(X;\,\Z/2)$ is onto for some
$x\in X$, then $\hat\theta(X)=0$.
Indeed, by the proof of Lemma 3.6(a), $\hat H^{2n}(K;\,\Z/2)=0$ for all
compact invariant $K\i\bar X$.
Hence from the Bockstein sequence, every element in $\hat H^{2n}(K)$ is
divisible by $2$.
Since it may be assumed finitely generated \cite{Bre}, it contains no
$2$-torsion.
So $\theta(X)$ restricts trivially over every $K$, thus $\hat\theta(X)=0$.
\endremark

\example{Example 2.7 (Sklyarenko's compactum)}
Fix any $n\ge 2$ and consider the ``local sphere'' $S^{n-1}_{(p)}$,
that is the mapping telescope of a direct sequence of degree $p$ maps
$S^{n-1}\to S^{n-1}$.
Let $X=X(n,p)$ be the one-point compactification of $S^{n-1}_{(p)}$, and let
$\infty$ denote the added point.
It is easy to see that $X$ is contractible and locally contractible.
It was noticed in \cite{Skl; Example 4.6} that $H^{n+1}(X,X\but\infty)$ is
non-zero.

Indeed, by the Milnor sequence this group is isomorphic to
$\lim^1 H^n(X,U_i)$, where $U_i$ is the union of the first $i$ mapping
cylinders in the mapping telescope.
Now each $G_i:=H^n(X,U_i)\simeq\Z$, and each support enlargement map
$G_{i+1}\to G_i$ is the multiplication by $p$, for the fundamental cocycle
$\omega_i$ with support in a point $p_i\in U_{i+1}\but U_i$ is cohomologous
with support in $X\but U_{i-1}$ to $p\omega_{i-1}$.
Hence $\lim^1\ne 0$.
\endexample

\remark{Remark}
Note that $X$ embeds in $\R^{2n}$ since there exists a pseudo-isotopy%
\footnote{We recall that a {\it pseudo-isotopy} of $\R^m$ is a homotopy
$H_t\:\R^m\to\R^m$ such that for each $t_0<1$, $H_{t_0}$ is a homeomorphism;
and the pseudo-isotopy is said to take $g\:X\to\R^m$ onto $f\:X\to\R^m$
if $H_0=\id$ and $H_1g=f$.}
taking the standard embedding $S^{n-1}\i S^{2n-1}$ onto the map
$S^{n-1}@>2>>S^{n-1}\i S^{2n-1}$.
This pseudo-isotopy can be constructed directly for $n=2$ and from the Zeeman
unknotting theorem (or also directly) for $n>2$.
\endremark

\example{Example 2.8} In fact, degree $p$ maps $S^{n-1}\to S^{n-1}$ in
Example 2.7 can be replaced by arbitrary maps of compact polyhedra
(or ANRs) $P^{n-1}_i\to P^{n-1}_{i+1}$ inducing a homomorphism with
nontrivial cokernel on $(n-1)$-homology modulo torsion.
Indeed, then $\dots\to H^{n-1}(P_2)\to H^{n-1}(P_1)$, which is isomorphic to
$\dots\to H^n(X,U_2)\to H^n(X,U_1)$, does not satisfy the Mittag-Leffler
condition, hence has a nontrivial $\lim^1$.
\endexample

Let us say that a compactum $X$ {\it movably embeds} in $\R^m$ if it
quasi-embeds in $\R^m$, and for each $\eps>0$ there exists a $\delta>0$
such that for each $\gamma>0$, every $\delta$-map $X\to\R^m$ is homotopic
through $\eps$-maps to a $\gamma$-map.
The following corollary to Theorem 2.2 ensures that movable
embeddability implies embeddability.

\proclaim{Corollary 2.9} An $n$-dimensional compactum $X$, $n>3$, embeds in
$\R^{2n}$ if and only if there exists a homotopy $h_t\:X\to\R^{2n}$,
$t\in [0,\infty)$, where each $h_t$ is a $\frac1t$-map.
\endproclaim

\demo{Proof} Let $K_i\i\bar X$ be the set of all $\{x,y\}$ with
$\dist(x,y)\le\frac1i$, then the mapping telescope
$\Tel(\dots\to K_2\to K_1)$ has the same cohomology as $\bar X$. \qed
\enddemo

Since the deleted product of a polyhedron is homotopy equivalent to
the simplicial deleted product, a compact $n$-polyhedron movably embeds in
$\R^{2n}$, $n>2$, iff it embeds there.
(Both this assertion and Corollary 2.9 remain true in the metastable range,
see \S4.)

\example{Example 2.10 (an AR that does not movably embed in $\R^{2n}$)}
Let $X=X(n,p)$ be the Sklyarenko compactum, with $\infty$ regarded as
a basepoint.
The wedge $Y:=X\vee X$ is a contractible locally contractible
$n$-dimensional compactum, embeddable into $\R^{2n}$.
We claim that $Y$ does not movably embed into $\R^{2n}$.

A sequence of $\frac1i$-maps $f_i\:Y\to\R^{2n}$ can be described as follows.
Each $f_i$ sends the basepoint $\infty$ to the center $c$ of a $2n$-ball
$B^{2n}$, and properly embeds each copy of the telescope $S^{n-1}_{(p)}$ into
$B^{2n}\but\{c\}\cong S^{2n-1}\x[0,\infty)$ in a level-preserving fashion.
The two copies of the $(n-1)$-sphere $S_i:=\Fr_X(U_i)$ are linked by $f_i$ in
$S^{2n-1}\x\{i\}$ with linking number $1$, and, writing $\pi$ for
the projection of $S^{2n-1}\x\R$ onto the first factor,
$\pi f_i(U_i\sqcup U_i)$ is contained in a regular neighborhood of
$\pi f_i(S_i\sqcup S_i)$ in $S^{2n-1}$, so that $f_i(S_{i-j}\sqcup S_{i-j})$
has linking number $p^{i-j}$ in $S^{2n-1}\x\{i-j\}$.
Next, in the levels between $i$ and $i+\frac12$ we set $f_i$ to be the track
of a generic unlinking homotopy $h_t\:S_i\sqcup S_i\to S^{2n-1}$, whose
only double point occurs at the level $i+\frac14$.
Finally, in the levels between $i+\frac12$ and $\infty$ we construct $f_i$
as in the remark following Example 2.7, so that $f_i(S_j\sqcup S_j)$
is the unlink in $S^{2n+1}\x\{j\}$ for each $j>i$.

Suppose that some $f_i$ is homotopic to a $\frac1{i+1}$-map through $1$-maps.
Let us consider $Z_j=(U_j\x X\cup X\x U_j)/T$.
So we have a homotopy $H_t\:X\to\R^{2n}$, whose only double points
$H_t(x)=H_t(y)$ occur for $\{x,y\}\in\bar Y\but Z_1$, from $f_i$ to
the map $H_1$ whose only double points occur for
$\{x,y\}\in\bar Y\but Z_{i+1}$.
Let $g_i\:X\emb\R^{2n+1}$ be an embedding projecting onto $f_i$.
Due to the existence of $H_t$,
$\bar g_i^*(\xi)=0\in H^{2n}(Z_{i+1},Z_1)$.
On the other hand, since $f_i$ has a unique double point at the level
$i+\frac14$, $\bar g_i^*(\xi)$ is represented by the fundamental cocycle
with support in a point $\{x,y\}\in\bar Y$ where $x$ and $y$ lie in the two
copies of $U_{i+1}\but U_i$.
It follows that the restriction
$H^{2n}(Z_{i+1},Z_1)\to H^{2n}(U_{i+1}\x U_{i+1},Z_1\cap U_{i+1}^2)\simeq\Z$
sends $\bar g_i^*(\xi)$ to $\pm p^{i-1}\x p^{i-1}\ne 0$, which is
a contradiction.
\endexample

\proclaim{Theorem 2.11} An $n$-dimensional compactum $X$, $n>3$, movably
embeds in $\R^{2n}$ if and only if $\hat\tta(X)=0$ and
$\lim^1 H^{2n-1}(K_i)=0$, where $K_i$ is a nested sequence of compact sets
in $\bar X=\bigcup K_i$.
\endproclaim

The proof is similar to that of Theorem 2.2 (see also \cite{AM1} and
\cite{M2}) and is left to the reader.

\example{Example 2.12 (an ANR, quasi-embeddable but non-embeddable in
$\R^{2n}$)} \linebreak Let $K$ be the $n$-skeleton of the $(2n+2)$-simplex
$\Delta^n*\Delta^{n+1}$, $n>2$.
Let $x$ be a point in the interior of the $n$-cell $\Delta^n$ of $K$ and
$B^n$ be a closed $n$-ball in the interior $U$ of another $n$-cell of $K$.
Let $Y$ be $K$ with $B^n$ replaced with a copy of the Sklyarenko compactum
$X=X(n,p)$, where $p$ is odd.
Let $Z$ be $Y$ with $x$ and the point $\infty\in X$ identified with each
other.
Then $Z$ is a locally contractible $n$-dimensional compactum.
We claim that $Z$ quasi-embeds but does not embed in $\R^{2n}$.

Indeed, by \cite{M5; Example 2.3}, there is a map $K\to\R^{2n}$ with one
transverse double point $c=d$, where $c$ lies in $U\but B^n$ and $d\ne x$
lies in the interior of $\Delta^n$.
By the preceding remark, one can convert this map into a map $Y\to\R^{2n}$
with one transverse double point $c=d$.
By general position one can pinch the latter to get a map $Z\to\R^{2n}$
with one transverse double point $c=d$.
Thus $\tta(Z)$ is the image of a generator
$\zeta\in H^{2n}(\bar Z,\bar Z\but\{c,d\})$.
Since $\tta(Z)$ has order two, it also equals the image of $p^i\zeta$ for
any given $i\in\N$.
But the image of $p^i\zeta$ can be represented by a cocycle with support
in $c\x(X\but U_i)$.
Thus $\tta(Z_i)=0$, where $Z_i$ is $Z$ with $X\but U_i$ identified to
a point.
Since $n>2$, each $Z_i$ embeds in $\R^{2n}$, whence $Z$ quasi-embeds in
$\R^{2n}$.

On the other hand, suppose that $g\:Z\emb\R^{2n}$ is an embedding.
Let $S^{n-1}_0$ denote $\Fr_ZX$ and $S^{n-1}_i$ denote $\Fr_X(U_i)$ for
$i>0$.
If $g|_{S^{n-1}_0\sqcup\partial\Delta^{n+1}}$ has linking number $l$,
then $g|_{S^{n-1}_i\sqcup\partial\Delta^{n+1}}$ has linking number $p^{-i}l$.
Since the latter must be an integer for every $i\in\N$, the only possibility
is $l=0$.
But this cannot be by \cite{M5; Lemma 2.5}.
\endexample

Theorem 2.2 has the following

\proclaim{Corollary 2.13} Any $\bmod2$ locally acyclic $n$-dimensional
compactum, $n>3$, quasi-embeddable in $\R^{2n-1}$, embeds in $\R^{2n}$.
\endproclaim

See discussion preceding Corollary 3.5 for a different proof of a special
case of Corollary 2.13.

\demo{Proof}
Let $X$ be the compactum, $\eta\:\tl X\to\bar X$ the cover, and $K_i\i\bar X$
the set of all unordered pairs $\{x,y\}$ with $\dist(x,y)\ge 1/i$.
Then a $1/i$-map $f_i\:X\to\R^{2n-1}$ yields a map
$\bar f_i\:K_i\to\RP^{2n-2}$, $\{x,y\}\mapsto\left<f(x)-f(y)\right>$,
and the restriction of $w_1(\eta)^{2n-1}$ to $H^{2n-1}(K_i)$ equals
$\bar f_i^*(\xi)$, which is zero.
Since $X$ is locally acyclic $\bmod 2$ and $K_i\i\Int K_{i+1}$,
the image of $H^{2n-2}(K_{i+1};\Z/2)$ in $H^{2n-2}(K_i;\Z/2)$ is finitely
generated \cite{Bre}, hence finite.
Then the inverse sequence of $H^{2n-2}(K_i;\Z/2)$ satisfies the
Mittag-Leffler condition, so its derived limit vanishes, hence
$w_1(\eta)^{2n-1}=0$ from the Milnor sequence.
But $\tta(X)=e(\eta)^{2n}$ is the Bockstein image of $w_1(\eta)^{2n-1}$
(cf\. \cite{M5}), so it vanishes. \qed
\enddemo

We note that among locally acyclic compacta are all homology manifolds.

\proclaim{Corollary 2.14} A $\bmod2$ locally acyclic inverse limit of
$n$-manifolds embeds in $\R^{2n}$, $n>3$, if either the manifolds are all
orientable, or $n$ is not a power of $2$.
\endproclaim

Duvall and Husch show that for $n=2^l$, an inverse limit of $3$-coverings
between connected sums of $\RP^n$'s, and the inverse limit of trivial
$2$-coverings between disjoint unions of $2^i$ copies of $\RP^n$'s (that is,
$\RP^n$ cross the Cantor set) do not embed into $\R^{2n}$ even up to shape
\cite{DH}.
These compacta are, of course, not locally acyclic.
However, the condition of local acyclicity is necessary in Corollary 2.14
by the following example.

\example{Example 2.15 (Akhmetiev's compactum)}
For each $n=2^l-1$ the $2$-cover $p\:S^n\to\RP^n$ composed with any
topological embedding $\RP^n\emb\R^{2n}$ is not approximable by topological
embeddings (for $n=1,3,7$ see \cite{M4; Examples 1 and 3} and for $l>3$
also \cite{DH; Example 4.8}).
It follows that the ``discrete mapping cylinder'' of $p$, that is,
$X:=S^n\x\N^*/_{(x,\infty)\sim (-x,\infty)}$, where
$\N^*=\N\cup\{\infty\}$
denotes the one-point compactification of $\N$, does not embed in $\R^{2n}$.
Note that $X$ is the inverse limit of the $n$-manifolds
$M_i:=S^n_1\sqcup\dots\sqcup S^n_i\sqcup\RP^n$ (which are parallelizable
for $n=1,3,7$) and the bonding maps $M_{i+1}\to M_i$ that restrict to
the identity on $M_i\i M_{i+1}$ and to the $2$-cover $S^n_{i+1}\to\RP^n$.
\endexample

This example answers \cite{M2; Problem in the end of \S3} in the
case $n=3,7$.
Partial results on the remaining case $n=2$ were obtained in \cite{CF}.

$\tta(X)$ can be computed using the following

\proclaim{Lemma 2.16 (geometric computation of $\tta$)} Let $X$ be an
$n$-dimensional compactum, quasi-embeddable in $\R^{2n}$.

(a) Let $f_i\:X\to\R^{2n}$ be an $\frac1i$-map for each $i\in\N$.
Let $K_i\i\bar X$ consist of all pairs $\{x,y\}$ with $\dist(x,y)\ge\frac1i$,
and pick some homotopy $H_i\:K_i\x I\to\RP^\infty$ between
$\bar f_i\:K_i\to\RP^{2n-1}\i\RP^\infty$ and the restriction of
$\bar f_{i+1}$.
Let $d(f_i,f_{i+1})\in H^{2n-1}(K_i)$ be the image of
$H_i^*(\Xi)\in H^{2n}(K_i\x I,K_i\x\partial I)$ under the Thom isomorphism,
where $\Xi\in H^{2n}(\RP^\infty,\RP^{2n-1})\simeq\Z$ is a fixed generator.
Then $\tta(X)\in\lim^1 H^{2n-1}(K_i)=\coker\tau$ is the coset of
$(d(f_i,f_{i+1}))\in\prod H^{2n-1}(K_i)$, where
$\tau\:\prod H^{2n-1}(K_i)\to\prod H^{2n-1}(K_i)$ is given by
$\tau(a_1,a_2,\dots)=(a_1-a_2|_{K_1},a_2-a_3|_{K_2},\dots)$.

(b) Suppose that $X$ is the inverse limit of $\{P_i;\,p^i_j\}$
so that each $p^\infty_i\:X\to P_i$ is a $\frac1i$-map and each $P_i$
admits an embedding $g_i$ into $\R^{2n}$.
Let $q^\infty_i\:K_i\to\tl P_i$ be the restriction of
$p^\infty_i\x p^\infty_i$ and let
$q^i_{i-1}\:\bar P_i\but\Sigma(p^i_{i-1})\to\bar P_{i-1}$ lift to the
restriction of $p^i_{i-1}\x p^i_{i-1}$, where $\Sigma(f)$ denotes
the set of all unordered pairs $\{x,y\}$ with $x\ne y$ and $f(x)=f(y)$.
Then $d(g_ip^\infty_i,g_{i-1}p^\infty_{i-1})=
(q^\infty_i)^*d(\bar g_i,\bar g_{i-1} q^i_{i-1})$, where $d(\phi,\psi)$ is
the first obstruction to homotopy of $\phi$ and $\psi$.%
\footnote{Given two maps
$\phi,\psi\:\bar P_i\but\Sigma(p^i_{i-1})\to\RP^{2n-1}$, classifying
the restriction of the $2$-cover $\tl P_i\to\bar P_i$ and coinciding on
the $(2n-2)$-skeleton, $d(\phi,\psi)\in
H^{2n-1}(\bar P_i\but\Sigma(p^i_{i-1});\pi_{2n-1}(\RP^{2n-1}))$
is defined to be to the class of the cocycle, assigning to a $(2n-1)$-cell
$\sigma$ the homotopy class of the map
$\phi\cup\psi\:\sigma\cup_{\partial\sigma}\sigma\to\RP^{2n-1}$.
If $\phi$ and $\psi$ are homotopic on the $(2n-2)$-skeleton, let $\psi'$
be homotopic to $\psi$ and coincide with $\phi$ on the $(2n-2)$-skeleton,
then $d(\phi,\psi):=d(\phi,\psi')$ is well-defined.}
\endproclaim

In particular, $\tta(X)=0$ if and only if the infinite system of equations
$$x_i-p(x_{i+1})=(q^\infty_i)^*d(\bar g_i,\bar g_{i-1} q^i_{i-1})$$
in the inverse sequence $\cdots@>p>>H^{2n-1}(K_2)@>p>>H^{2n-1}(K_1)$ is
soluble.

\demo{Proof} Part (b) is clear.
Let $T$ be the telescope of the sequence of inclusions
$K_1\i K_2\i\cdots$, and let $p\:T\to\bar X$ be the projection and
$H\:(T,\bigsqcup K_i\x i)\to(\RP^\infty,\RP^{2n})$ be obtained by combining
the $H_i$'s.
Then $p^*\:H^{2n}(\bar X)\to H^{2n}(T)$ is an isomorphism, which sends
$\tta(X)$ to $H^*(\xi)$, where $\xi\in H^{2n}(\RP^\infty)$ is the generator.
On the other hand, $H^*(\xi)$ is the image of $H^*(\Xi)$ under the
homomorphism $H^{2n}(T,\bigsqcup K_i\x i)\to H^{2n}(T)$.
Hence $H^*(\xi)=0$ if and only if $H^*(\Xi)$ lies in the image of
$\delta^*\:H^{2n-1}(\bigsqcup K_i\x i)\to H^{2n}(T,\bigsqcup K_i\x i)$.
The latter can be identified with $f$, whence the coset of $H^*(\Xi)$ in
$\coker\delta^*$ gets identified with $\tta(X)\in\coker\tau$. \qed
\enddemo

\remark{Remark} The statement of Lemma 2.16 may look cumbersome, but it
formalizes a simple idea.
If $X$ is a polyhedron, $\theta(X)\in H^{2n}(\bar X)$ can be represented by
the double point set $\Sigma(f)\i\bar X$ of a generic PL map $X\to\R^{2n}$
(see \cite{M5}).
To be precise about geometrically representing cohomology classes, the
representatives must be thought of as mock bundles (with codimension two
singularities, which however do not occur in the present context, for
a good reason \cite{M5; \S4}).
If $X$ is a compactum that is sufficiently ``tame'' to be considered as
a polyhedron with singularities (e.g\. if it is a $Z$-compactification of
an infinite polyhedron), and it admits a ``generic'' map $f$ into $\R^{2n}$,
whose double points do not directly involve the singularities, one can
apply the same method to $f$.
Another option, which works for arbitrary $X$, is to apply this method
to a level-preserving generic map of the mapping telescope of an inverse
sequence of nerves of $X$ into $\R^{2n}\x [0,\infty)$.
Writing this out produces the statement of Lemma 2.16, where the difference
$d(\bar g_i,\bar g_{i-1} q^i_{i-1})$ is now interpreted geometrically,
as represented by the double point set of a generic homotopy between
$g_{i-1} q^i_{i-1}$ and $g_i$, restricted to $\bar P_i\but\Sigma(p^i_{i-1})$.
\endremark

\example{Example 2.17 (the Ljubljana compactum)}
Let $T$ be the cone over the $(n-1)$-skeleton of the $2n$-simplex
(e.g. $T$=triod when $n=1$), and let $\N^*=\N\cup\{\infty\}$ denote
the one-point compactification of $\N$.
Let us show that the contractible $n$-dimensional compactum
$T\x\N^*\cup c\x [0,\infty]$ does not embed in $\R^{2n}$, reproving the
result of \cite{RSS} (see also \cite{RS}, \cite{KR}).
The same argument works if $T$ is an $n$-dimensional product of
the cones over the $n_i$-skeleta of $(2n+2)$-simplices (e.g. the product of
$n$ triods).

It suffices to show that $T\x\N^*$ does not embed in $\R^{2n}$.
Now $\overline{T\x\N^*}$ contains the union $K$ of the increasing sequence of
compacta $$K_i\ \ =\ \ \bar T_s\x\Delta_{\N^*}\ \cup\ \tl T_s\x\bar\N_i^*\
\cup\ T^2\x(\bar\N^*\but\bar\N_i^*),$$ where $\N_i^*=\{i,i+1,\dots,\infty\}$
and $s$ stands for the simplicial deleted product (see \cite{M5}).
By functoriality of the Euler class it suffices to show that
$e(\eta)^{2n}\ne 0$, where $\eta$ denotes the restriction of the $2$-cover
$\widetilde{T\x\N^*}\to\overline{T\x\N^*}$ over $K$ (see \cite{M5}).
Since $\N^*$ is the inverse limit of the finite sets $\{1,2,\dots,j\}$,
each $K_i$ is the inverse limit of the compact polyhedra obtained from $K_i$
by replacing each occurrence of $\N^*$ (resp\. $\N^*_i$) with
$\{1,2,\dots,j\}$ (resp\. $\{i,i+1,\dots,j\}$).
Since the $2$-cover $\tl T_s\to\bar T_s$ is homeomorphic to the $2$-cover
$S^{2n-1}\to\RP^{2n-1}$ (see \cite{M5}), it follows that $H^{2n-1}(K_i)$
can be identified with the subgroup $G_i$ of
$\left(\prod_{j\ge 1}2\Z\right)\x\prod_{j\ge i}\prod_{k>j}\Z$
consisting of all integer sequences $(n_{jk})$, where $k\ge j$ and
either $j\ge i$ or $k=j$, with each $n_{jj}\in 2\Z$ and eventually
stabilizing as $k\to\infty$ in the sense that $\exists k_0$
$\forall j\ge i$ $\forall k\ge\max(j,k_0)$ $n_{jk}=n_{\min(j,k_0),k_0}$.
(Thus abstractly $G\simeq\bigoplus_{l\ge 1}\Z$.)
We have the commutative diagram
$$\CD\dots@>p_2>>G_2@>p_1>>G_1\\
@.@VVV@VVV\\
\dots@>>>\left(\prod\limits_{j\ge 1}2\Z\right)\x\prod\limits_{j\ge 2}
\prod\limits_{k>j}\Z@>>>
\left(\prod\limits_{j\ge 1}2\Z\right)\x\prod\limits_{j\ge 1}
\prod\limits_{k>j}\Z,\!
\endCD$$
and $H^{2n-1}(K)$ is isomorphic to the inverse limit of the upper row.

Similarly, each $G'_i:=H^{2n-1}(K_i;\,\Z/2)$ can be identified with
the subgroup of all eventually stabilizing sequences in
$\left(\prod_{j\ge 1}2\Z/4\right)\x\prod_{j\ge i}\prod_{k>j}\Z/2$.
It is easy to see that
$w_1(\eta)^{2n-1}\in H^{2n-1}(K;\,\Z/2)\simeq\lim G'_i$ has all
$n_{jj}=2+4\Z\in 2\Z/4$ and all $n_{jk}=0$ whenever $k>j$.

We will now compute $e(\eta)^{2n}$ using Lemma 2.18 below.%
\footnote{A computation of $e(\eta)^{2n}$ based on Lemma 2.16 may be shorter,
if the proof of Lemma 2.18 is taken into account.
The arguments in \cite{RSS}, \cite{RS}, \cite{KR} are also not hard; yet
another option is to induce $\tta(T\x\N^*)$ from $\tta(T\x I)$.
The point here is to illustrate Lemma 2.18, which may find other applications
in the context of Problem 1.1.}
First we find a lift of $w_1(\eta)$ into some
$(n_{jk}^1,n_{jk}^2,\dots)\in\prod G_i$.
We may take all $n_{jj}^i=2\in 2\Z$ and all $n_{jk}^i=2\in\Z$ when $k>j$.
(We may not take all $n_{jj}^i=2$ and all $n_{jk}^i=0$ when $k<j$, since
such a sequence does not eventually stabilize.)
Now let $m^i=n^i-p_i(n^{i+1})$, i.e\. $m^i_{jk}=n^i_{jk}-n^{i+1}_{jk}$
when $n^{i+1}_{jk}$ is defined, and $m^i_{jk}=n^i_{jk}$ otherwise.
Then $m^i_{jk}=0$ when $j>i$ or $k=j$ and $m^i_{ik}=2$ for all $k>i$.
By Lemma 2.18, $e(\eta)^{2n}\in\lim^1 G_i$ is represented by
$(m^1/2,m^2/2,\dots)$.
It is not hard to see that this element is nontrivial, that is,
the infinite system of equations $x_i-p_i(x_{i+1})=m^i/2$ has no solution.
The system $x_i-p_i(x_{i+1})=m^i$ admits the solution $x_i=n^i$, in
accordance with $2e(\eta)^{2n}=0$.
\endexample

\proclaim{Lemma 2.18 (intrinsic computation of $\tta$)} Let $\eta$ be
a $2$-cover of the union $K$ of a sequence $K_1\i K_2\i\cdots$, and suppose
that $H^{2n}(K_i)=0=H^{2n-2}(K_i;\,\Z/2)$ for each $i$.
Then the element $e(\eta)^{2n}\in\lim^1 H^{2n-1}(K_i)$
is the image of $w_1(\eta)^{2n-1}\in\lim H^{2n-1}(K_i;\,\Z/2)$ under
the connecting homomorphism from the six-term exact sequence corresponding
to the short exact sequence of inverse sequences
$$0@>>>H^{2n-1}(K_i)@>2>>H^{2n-1}(K_i)@>>>H^{2n-1}(K_i;\,\Z/2)@>>>0.$$
\endproclaim

Note that the nature of $\eta$ is irrelevant for the conclusion.

\demo{Proof} Let $T$ be the telescope of the sequence of inclusions
$K_1\i K_2\i\cdots$, and let $p\:T\to K$ be the projection.
Then $p^*\:H^m(K)\to H^m(T)$ is an isomorphism, and $\lim H^m(K_i)$ is
the image of $H^m(T)$ in $H^m(\bigsqcup K_i\x i)\simeq\prod H^m(K_i)$.
The image of $\prod H^{m-1}(K_i)\simeq H^m(T,\bigsqcup K_i\x i)$
in $H^m(T)$ is isomorphic to $\lim^1 H^{m-1}(K_i)$ (see proof of Lemma 2.16).

On the other hand, $\lim H^{2n-1}(K_i)$ and $\lim^1 H^{2n-1}(K_i)$ are
the cohomology groups of the $1$-dimensional cochain complex
$\prod H^{2n-1}(K_i)@>\tau>>\prod H^{2n-1}(K_i)$ with coboundary
$\tau(a_1,a_2,\dots)=(a_1-a_2|_{K_1},a_2-a_3|_{K_2},\dots)$.
This $\tau$ can be identified with
$\delta^*\:H^{2n-1}(\bigsqcup K_i\x i)\to H^{2n}(T,\bigsqcup K_i\x i)$.
So what we need to show is that
$e(\eta)^{2n}\in\im\left[H^{2n}(T,\bigsqcup K_i\x i)\to H^{2n}(T)\right]$
lifts to an $y\in H^{2n}(T,\bigsqcup K_i\x i)$ such that
$2y=\delta^*(z)$, where $z\bmod 2\in H^{2n-1}(\bigsqcup K_i\x i;\,\Z/2)$
coincides with $w_1(\eta)^{2n-1}\in\im\left[H^{2n-1}(T;\,\Z/2)\to
H^{2n-1}(\bigsqcup K_i\x i;\,\Z/2)\right]$.

Now $e(\eta)^{2n}\in H^{2n}(K)$ is the Bockstein image of
$w_1(\eta)^{2n-1}\in H^{2n-1}(K;\,\Z/2)$ (see \cite{M5}).
Hence $e(p^*\eta)^{2n}\in H^{2n}(T)$ is the Bockstein image of
$w_1(p^*\eta)^{2n-1}\in H^{2n-1}(T;\,\Z/2)$.
In other words, $e(p^*\eta)^{2n}=[Y]$ where $2Y=\delta Z$ where
$[Z\bmod 2]=w_1(p^*\eta)^{2n-1}$.

Since $H^{2n}(T,\bigsqcup K_i\x i)\to H^{2n}(T)$ is onto, we may assume that
$Y$ vanishes on $L:=\bigsqcup K_i\x i$; set $y=[Y]$.
Then $\delta(Z|_L)=\delta Z|_L=Y|_L=0$, i.e\. $Z|_L$ is a cocycle;
set $z=[Z|_L]$.
Finally, $\delta^*(z)=[X]$ where $X=\delta W$ where $W|_L=Z|_L$.
We may take $W=Z$, thus $\delta^*(z)=[\delta Z]=2[Y]=2y$. \qed
\enddemo

\head 3. Proofs \endhead

\proclaim{Proposition 3.1} Let $X$ be the limit of an inverse sequence
$\{X_i;p^i_j\}$ of compact $n$-polyhedra and PL maps, $n>1$.
If $\hat\tta(X)=0$, then for each $k$ and $\eps>0$ there exists an $l>k$
such that the composition of $p^l_k$ and an $\eps$-map $X_k\to Y$, factors
through a compact $n$-polyhedron $K$ with $\tta(K)=0$.
\endproclaim

\demo{Proof}
If $\phi$ is a map to $X_k$, let $\Sigma_\eps(\phi)$ denote
$\{\{x,y\}\mid\dist(\phi(x),\phi(y))<\eps\}$.
By naturality of the Euler class, the image of $\hat\tta(X)$ in
$H^{2n}(\bar X\but\Sigma_\eps(p^\infty_k))
=\dirlim H^{2n}(\bar X_i\but\Sigma_\eps(p^i_k))$
is the thread consisting of the restriction images of $\tta(X_i)$, $i\ge k$.
Since $\hat\tta(X)=0$, there exists an $l>k$ such that the restriction image
of $\tta(X_l)$ in $H^{2n}(\bar X_l\but\Sigma_\eps(p^l_k))$ is zero.
In other words, $\tta(X_l)$ can be represented by a cocycle with support
in $\Sigma_\eps(p^l_k)$.
Since this cocycle is top-dimensional, its support can be chosen in a
finite set $S/T\i\Sigma_\eps(p^l_k)$, disjoint from the $(2n-1)$-skeleton
of some triangulation of $\bar X_l$.

After a small perturbation of $S/T$ we may assume that its $2$-cover
$S\i\tl X_l$ projects injectively onto the first factor of $X_l\x X_l$
and the set $R:=(p^l_k\x p^l_k)(S)\but\Delta_{X_k}$ projects injectively
onto the first factor of $X_k\x X_k$.
For each $(x,y)\in R$ we have $\dist(x,y)<\eps$.
Viewing $S$ and $R$ (or rather $S\cup\Delta_{X_l}$ and $R\cup\Delta_{X_k}$)
as equivalence relations on $X_l$ and $X_k$, we have that the composition
of $p^l_k$ and the quotient $\eps$-map $\pi\:X_k\to X_k/R$ factors through
the quotient $K:=X_l/S$.

The quotient map $X_l\to K$ induces a surjection
$q\:\bar X_l\but S\to\bar K$, which identifies $\{x,a\}$ with $\{x,b\}$
whenever $\{a,b\}\in S$.
The cone of $q$ collapses onto an $(n+1)$-polyhedron, and therefore
$q$ induces an isomorphism on $2n$-cohomology as long as $n\ge 2$.
The image of $\tta(K)$ under the latter coincides by naturality of
the Euler class with the image of $\tta(X_l)$ under the restriction
$H^{2n}(\bar X_l)\to H^{2n}(\bar X_l\but S/T)$.
Since the latter is zero by the above, $\tta(K)=0$. \qed
\enddemo

\remark{Remark} One may wonder whether it is possible to do without
the $\eps$-map in Proposition 3.1.
It would certainly be superfluous if we knew that
(i) $H^{2n}(\bar X\but\Sigma(p^\infty_k))$ is the direct limit of
$H^{2n}(\bar X_i\but\Sigma(p^i_k))$, and moreover,
(ii) the image of $\tta(X)$ in $H^{2n}(\bar X\but\Sigma(p^\infty_k))$
(rather than just the image of $\hat\tta(X)$ in
$\hat H^{2n}(\bar X\but\Sigma(p^\infty_k))$) is zero.

As for (i), the letter X (viewed as a finite graph, namely the cone over $4$
points) is the limit of an inverse sequence of the letters H, with bonding
maps shrinking the horizontal segment onto an increasingly small subsegment
in the center and fixing the four boundary points $\partial$H.
Similarly, the double $X$ of the letter X (that is, the suspension over $4$
points) is the inverse limit of the doubles
$X_i:=$H$\cup_{\partial\text{H}}$H of the letters H and the doubled bonding
maps $p^i_j$.
Consider a nonzero element
$\alpha\in H^2(\bar X_1)$ represented by a cocycle with support in $\{l,r\}$,
where $l$ (resp\. $r$) is the left (resp\. right) lower endpoint of the
letter H, viewed as contained in $\partial$H$\i X_1$.
The image of $\alpha$ in $H^2(\bar X\but\Sigma(p^\infty_1))$ is zero, but its
images in $H^2(\bar X_i\but\Sigma(p^i_1))$ are all nonzero.
Hence the assertion (i) does not hold.
Moreover, it is clear that no bonding map $p^i_1$ factors through
a map $p\:K\to X_1$ such that $\alpha$ trivializes in
$H^2(\bar K\but\Sigma(p))$.

Concerning (ii), note that $\dirlim$ does not always commute with
$\invlim$.
For instance,
$$\dirlim(\prod_{k\ge 1}\Z@>>>\prod_{k\ge 2}\Z@>>>\dots)\simeq
\left.\prod\Z\right/\bigoplus\Z\ne 0,$$ whereas
$\dirlim(\prod_{k=1}^i\Z@>>>\prod_{k=2}^i\Z@>>>\dots)=0$ for each $i\in\N$.
\endremark
\medskip

\proclaim{Criterion 3.2} A compactum $X$ embeds in $\R^m$ iff for each
$i\in\N$ there exist a $\frac1i$-map $f_i\:X\to\R^m$ and a pseudo-isotopy
$H^i_t\:\R^m\to\R^m$ taking $f_{i+1}$ onto $f_i$.
\endproclaim

Antecedents of this criterion and its proof include those in \cite{I} and
\cite{SS}.

\remark{Remark}
The following was proved in \cite{SS} using Shtanko's theory of embedding
dimension \cite{Sht}.
If $m-n\ge 3$ and an $n$-dimensional compactum $X$ embeds in $\R^m$, then
$X$ can be decomposed into an inverse sequence $\{X_i;p^i_j\}$ of compact
$n$-polyhedra and PL maps such that there exist PL embeddings
$f_i\:X_i\to\R^m$ and PL pseudo-isotopies $H^i_t\:\R^m\to\R^m$ taking
$f_{i+1}$ onto $f_ip_i^{i+1}$.
\endremark

\demo{Proof} Consider the maps $g_1=f_1$, $\,g_2=H^1_{\eps_1}f_2$,
$\,g_3=H^2_{\eps_2}H^1_{\eps_1}f_3$, etc., where the $\eps_i$'s are chosen
so small that each $g_{i+1}$ is $\delta_i$-close to $g_i$, where
$\delta_i=\min\{\frac{\gamma_j}{3^{i+1-j}}\mid 1\le j\le i\}$
and $\gamma_j>0$ is such that $\dist(g_j(x),g_j(y))\ge 2\gamma_j$
whenever $\dist(x,y)\ge\frac1j$.
Since $\delta_1+\delta_2+\dots$ has a convergent majorant
$\gamma_1(\frac13+\frac19+\dots)$, the $g_i$'s uniformly converge to
a map $g\:X\to\R^m$.
Since $\delta_i+\delta_{i+1}+\dots\le\frac{\gamma_i}2$, for any points
$x,y\in X$ with $\dist(x,y)\ge\frac1i$, the triangle axiom implies
$\dist(g(x),g(y))\ge 2\gamma_i-\frac{\gamma_i}2-\frac{\gamma_i}2>0$.
Thus $g$ is an embedding. \qed
\enddemo

The content of Criterion 3.2 is illustrated by the following simple
observation.

\proclaim{Proposition} If $P$ is a compact $n$-polyhedron, for every PL map
$f\:P\to\R^{2n+1}$ there exists a PL embedding $g\:P\to\R^{2n+1}$ and
a PL pseudo-isotopy of $\R^{2n+1}$ taking $g$ onto $f$.
\endproclaim

\remark{Remark} The following addendum to this proposition was proved in
\cite{M1}: if $n>1$, then for each $\eps>0$ there exists a $\delta>0$
such that any PL embedding, $\delta$-close to $f$, can be taken onto $f$ by
a PL $\eps$-pseudo-isotopy.
This is not the case for $n=1$ \cite{M1}.
The proposition and the addendum hold even if $\R^{2n+1}$ is replaced with
any $\R^m$, where $m-n\ge 3$ \cite{M1}.
On the other hand, the proposition and the addendum remain true if
one drops all occurrences of ``PL'', provided that $n>1$ \cite{M1}.
When $n=1$, this is not so for the addendum (it can be shown, using the first
author's multivariable version of the Conway polynomial, that the Bing sling
precomposed with the trivial $2$-cover of $S^1$ cannot be obtained by
a pseudo-isotopy of any embedding), and whether this is so for the
proposition depends on a certain conjecture in classical (tame) link theory.
\endremark

\demo{Proof} In some triangulations where $f$ is simplicial, let $Q$
be the $2n$-skeleton of the dual cell complex $\Cal C$ of $\R^{2n+1}$, and
let $R$ denote the $(n-1)$-polyhedron $f^{-1}(Q)$.
By general position, there exists a PL embedding $g|_R\:R\to Q$ and a PL
pseudo-isotopy $H_t|_Q\:Q\to Q$ such that $H_0=\id$ and $H_1(g|_R)=f|_R$.
Extend $g|_R$ by general position to an arbitrary PL embedding
$g\:P\emb\R^{2n+1}$ such that $g^{-1}(v_i^*)=f^{-1}(v_i^*)$ for each
$(2n+1)$-cell $v_i^*$ of $\Cal C$ dual to a vertex $v_i$ of
the triangulation of $\R^{n+1}$.
Viewing $H_t$ as $H\:Q\x I\to Q\x I$, extend it conewise to each cone
$(v_i,1)*(\partial v_i^*\x I)$, whose ``inner'' boundary is
$D_i:=(v_i,1)*(\partial v_i^*\x 0)$, and then conewise to each cone
$(v_i,0)*D_i$.
Then $H_1g=f$. \qed
\enddemo

\definition{Notation}
If $f\:X\to Y$ is a map, let $\Sigma(f)\i\bar X$ denote the set of all
unordered pairs $\{x,y\}$ of distinct points of $X$ such that $f(x)=f(y)$.
If an $(n+1)$-dimensional $\sigma$-compactum $T$ is naturally equipped with
a proper map $\del\:T\to [0,\infty)$ (its choice will be clear from
the context in each case), we define $\bar T^\del=\Sigma(\del)\i\bar T$,
its double cover $\tl T^\del\i\tl T$, and
$\tta_\del(T)=e(\tl T^\del\to\bar T^\del)^{2n}$.
\enddefinition

\proclaim{Lemma 3.3} Let $T$ be the mapping telescope of an inverse sequence
$\{X_i;p^i_j\}$ of compact $n$-polyhedra and PL maps, $n>3$.
If $\tta_\del(T)=0$, then for each $i$ there exists a PL embedding
$g_i\:X_i\emb\R^{2n}$ and a PL pseudo-isotopy taking $g_{i+1}$ onto $g_i$.
\endproclaim

\demo{Proof} We have $G^*(\xi)=0$ for some $G\:\bar T^\del\to\RP^{2n+1}$
classifying $\eta\:\tl T^\del\to\bar T^\del$, where
$\xi\in H^{2n}(\RP^{2n+1})$ is the generator.
Let $U$ be some triangulation of the $(2n+1)$-polyhedron $\bar T^\del$.
We may assume that $G(U^{(2n-1)})\i\RP^{2n-1}$.
Then $G^*(\xi)$ is by definition the obstruction to the existence of a map
$U^{(2n)}\to\RP^{2n-1}$, coinciding with $G$ on $U^{(2n-2)}$.
By the obstruction theory, such a map exists; if $n>1$, it still classifies
$\eta$.
The obstruction to the existence of a map $U\to\RP^{2n-1}$, coinciding on
$U^{(2n-1)}$ with the map just constructed, lies in
$H^{2n+1}(U;\pi_{2n}(\RP^{2n-1}))$.
(The coefficient sheaf is constant since its stalks have order two.)

This group can be computed from the Milnor exact sequence
$$0\to\derlim H^{2n}(U_i;\Z/2)\to H^{2n+1}(U;\Z/2)\to\invlim
H^{2n+1}(U_i;\Z/2)\to 0,$$
where $U_i=\del^{-1}([0,i])$.
The inverse sequence on the left hand side consists of finite groups,
hence satisfies the Mittag-Leffler condition and therefore has trivial
derived limit.
On the other hand, each $U_n$ collapses onto the $(2n)$-polyhedron
$$\bar X_0\cup\bigcup_{i=1}^n N_i\cup\Cyl\left[(p^i_{i-1}\x p^i_{i-1})
|_{\Fr N_i}\right],$$
where each $N_i$ is the second derived neighborhood of
``$\Sigma(p^i_{i-1})\cup$infinity'' in some triangulations of $\bar X_i$
where $p^i_{i-1}$ are simplicial for $i\le n$.
Thus $H^{2n+1}(U;\Z/2)=0$.
By the obstruction theory there exists a map $F\:U\to\RP^{2n-1}$ classifying
$\eta$.
By the covering theory, $F$ lifts to a map $\tl T^\del\to S^{2n-1}$; it is
necessarily equivariant if $\tl T^\del$ is connected, and otherwise can be
easily made equivariant.

Now $\tl T^\del$ contains the mapping telescope of the inverse sequence of
polyhedra $\tl X_i$ and partial maps
$q^{i+1}_i\:\tl X_{i+1}\supset\tl X_{i+1}\but\Sigma(p^{i+1}_i)\to\tl X_i$
obtained by restricting $p^{i+1}_i\x p^{i+1}_i$.
Thus there is a sequence of equivariant maps $F_i\:\tl X_i\to S^{2n-1}$ such
that each $F_{i+1}$ restricted to $\tl X_{i+1}\but\Sigma(p^{i+1}_i)$ is
equivariantly homotopic to $F_i q^{i+1}_i$.
By the Haefliger--Weber criterion there exist embeddings
$g_i\:X_i\emb\R^{2n}$ such that each $\tl g_i$ is equivariantly homotopic to
$F_i$.
It follows that each $\tl g_{i+1}$ restricted to
$\tl X_{i+1}\but\Sigma(p^{i+1}_i)$ is equivariantly homotopic to
$\widetilde{g_ip^{i+1}_i}$.
Then by the rel$\partial$ version of the Haefliger--Weber criterion,
a computation of the deleted product of a mapping cylinder and some
PL topology \cite{Sk2; Pseudo-isotopy Theorem 5.5}+\cite{M1; Theorem 1.12}
there exists a PL pseudo-isotopy taking $g_{i+1}$ onto $g_i$.
(The proof of Theorem 5.5 in \cite{Sk2} contains a minor gap: it is not
proved there that a PL pseudo-concordance can be split into a
PL pseudo-isotopy and a genuine PL concordance --- a fact less
obvious than its DIFF counterpart; luckily, it was verified in \cite{M1}.)
\qed
\enddemo

\proclaim{Proposition 3.4} Let $X$ be the limit of an inverse sequence
$\{X_i;p^i_j\}$ of compact $n$-polyhedra and PL maps.
Let $T_1$ be its mapping telescope, and if $h\:\N\to\N$ is an increasing
sequence, let $T_h$ denote the mapping telescope of
$\{X_{h(i)};p^{h(i)}_{h(j)}\}$ and let $s^h_1\:T_h\to T_1$ be the fiberwise
map obtained by combining the $p^{h(i)}_i$'s.

(a) If $\tta(X)=0$, then for any proper map $\eps\:\N\to (0,1]$ there exists
an increasing sequence $l\:\N\to\N$ such that after amending each
$p^{i+1}_i$ by a self-map of $X_i$, $\eps(i)$-close to the identity,
$\tta_\del(T_l)$ can be represented by a cocycle with support in
$\Sigma(s^l_1)$.

(b) If $n>3$ and $\tta_\del(T_l)$ can be represented by a cocycle with
support in $\Sigma(s^l_1)$, then $s^l_1$ fiberwise factors through
the mapping telescope $K$ of an inverse sequence $\{K_i;q^i_j\}$ of compact
$n$-polyhedra and PL maps such that $\tta_\del(K)=0$.

(c) There exists a proper map $\eps\:\N\to (0,1]$ such that the inverse limit
of any $\{K_i;q^i_j\}$ satisfying the above with this $\eps$ is homeomorphic
to $X$.
\endproclaim

\demo{Proof. (a)}
Since each $p^{i+1}_i$ is uniformly continuous, without loss of generality
$\dist(x,y)<\eps(i+1)$ for $x,y\in X_{i+1}$ implies
$\dist(p^{i+1}_i(x),p^{i+1}_i(y))<\eps(i)$.
If $\phi$ is a map to $T_1$, let $\Sigma_\eps(\phi)$ denote the set of
all unordered pairs $\{x,y\}$ such that $\del(\phi(x))=\del(\phi(y))=:t$
and $\dist(\phi(x),\phi(y))<\eps(\lceil t\rceil)$.

{\bf Step 1: $\exists l$ such that $\tta_\del(T_l)$ can be represented by
a cocycle with support in $\Sigma_\eps(s^l_1)$.}
Write $T_\infty=X\x[0,\infty)$, and let $s^\infty_1\:T_\infty\to T_1$ be
the projection.
By naturality of the Euler class, the image of $\tta_\del(T_\infty)$ in
$H^{2n}(\bar T_\infty^\del\but\Sigma_\eps(s^\infty_1))=
\dirlim H^{2n}(\bar T_h^\del\but\Sigma_\eps(s^h_1))$
(the direct limit over all increasing sequences $h\:\N\to\N$)
is the thread consisting of the restriction images of $\tta_\del(T_h)$.
Since $\tta(X)=0$, we have $\tta_\del(T_\infty)=0$, hence there exists an
increasing sequence $l\:\N\to\N$ such that the restriction image of
$\tta_\del(T_l)$ in $H^{2n}(\bar T_l^\del\but\Sigma_\eps(s^l_1))$ is zero.
In other words, $\tta_\del(T_l)$ can be represented by a cocycle with support
in $\Sigma_\eps(s^l_1)$.

Since this cocycle is of codimension one, its support can be chosen in
a $1$-polyhedron $S/T\i\Sigma_\eps(s^l_1)$, where $S\cup\Delta_{T_l}$ is
a subpolyhedron of $T_l\x T_l$, and $S\i\tl T_l^\del$ is disjoint from
the $(2n-1)$-skeleton
$\bigcup_i(\alpha X_{l(i)}\x\alpha X_{l(i)})^{(2n-2)}\x (i-1,i]$ for some
triangulations $\alpha X_i$ of $X_i$ such that each $p^{i+1}_i$ is simplicial
in $\alpha X_{i+1}$ and some subdivision of $\alpha X_i$.
We may assume that $\eps(i)$ is less than the minimal distance between
disjoint simplices of $\alpha X_i$.
Let $X_l$ denote $\bigcup X_{l(i)}\x\{i\}\i T_l$ and let $1^+\:\N\to\N$ be
given by $1^+(i)=i+1$.
Without loss of generality $l(i-1)>i$ for each $i$, so
$s^l_{1^+}\:T_l\to T_{1^+}$ is defined.

{\bf Step 2: $S/T\cap\bar X_l^\del\i\Sigma(s^l_{1^+})$ after amending each
$p^{i+1}_i$ by a self-map of $X_i$, $\eps(i)$-close to the identity.}
Fix some $i\in\N$.
Without loss of generality, $A_i:=S\cap\tl X_{l(i-1)}$ is a finite set.
After a small perturbation (of $S/T$) we may assume that
$B_i:=p^{l(i-1)}_i\x p^{l(i-1)}_i(A_i)$ projects injectively to
the first factor of $X_i\x X_i$.
Since $A_i\i\Sigma_{\eps(i)}(p^{l(i-1)}_i)$, for each $(x,y)\in B_i$
we have $\dist(x,y)<\eps(i)$.
Since $\eps(i)<\mesh(\alpha X_i)$, the projection $B_i\i X_i\x X_i\to X_i$
extends to a map $\lambda_i\:\Cyl(B_i\to B_i/T)\to X_i$, which sends
$B_i/T$ to the set $\alpha^*X_i^{(1)}\cap\alpha X_i^{(n-1)}$ of
the barycenters of all $(n-1)$-simplices of $\alpha X_i$, and embeds
the complement of $B_i/T$ into the complement of $\alpha X_i^{(n-1)}$.
Then the image of $\lambda_i$ is a disjoint union of $1$-polyhedra $C_{ij}$
(specifically, cones over finite sets), each of diameter $<\eps(i)$,
collapsible onto their intersections with $\alpha X_i^{(n-1)}$.
Hence each $C_{ij}$ is {\it shrinkable}, that is, if $N$ is its second
derived neighborhood in $\alpha X_i$, then $(N,\Fr N)$ is
PL homeomorphic to $(N/C_{ij},\Fr N)$.
Viewing $C_i:=\bigcup C_{ij}\x C_{ij}$ (or rather $C_i\cup\Delta_{X_i}$) as
an equivalence relation on $X_i$, consider the quotient $\eps(i)$-map
$\pi_i\:X_i\to X_i/C_i$.
Since each $C_{ij}$ is shrinkable, $\pi_i$ is arbitrarily closely
approximable by a homeomorphism $h_i\:X_i\to X_i/C_i$, and since
each $C_{ij}$ has diameter $<\eps(i)$, $g_i:=h_i^{-1}\pi_i$ is
a self-map of $X_i$, $\eps(i)$-close to the identity.
Writing $\acute p^{i+1}_i$ for $g_ip^{i+1}_i$, we have
$A_i\i\Sigma(\acute p^{i+1}_i)$ since $B_i\i C_i$.
Writing $\acute p^j_i$ for $\acute p^{i+1}_i\dots\acute p^j_{j-1}$, we also
have $A_i\i\Sigma(\acute p^{l(i-1)}_i)$ since $A_i$ may be assumed disjoint
from $(p^{l(i-1)}_j\x p^{l(i-1)}_j)^{-1}(C_j)$ for each $j$ with
$l(i-1)>j>i$.

Let $\acute T_l$ be the mapping telescope of
$\{X_{l(i)},\acute p^{l(j)}_{l(i)}\}$, and define
$\acute s^l_h\:\acute T_l\to\acute T_h$ by combining
the $\acute p^{l(i)}_{h(i)}$'s.
We have $S/T\cap\bar X_l^\del\i\Sigma(\acute s^l_{1^+})$.
To keep track of the entire $S$, we need to devise a map $T_l\to\acute T_l$.
The bonding maps $\acute p^j_k$ admit another definition in the case
where $j$ is known as a function of $k$, in particular when $k=l(i)$
and $j=l(i+1)$ for some $i$.
Let $C_k^{[j]}=C_k\cup p^{k+1}_k(C_{k+1})\cup\dots\cup p^{j-1}_k(C_{j-1})$.
We may assume that $C_k^{[j]}$ is collapsible onto finitely many points,
and is disjoint from $X_k^{(n-1)}$, except possibly at some of these points.
(To achieve this, we need to make sure e.g\. that $\bar C_{k+h}$ is disjoint
from $\Sigma(p^{k+h}_k)$ for each $h<j-k$, which requires knowing $j=j(k)$
in advance.)
Then we may use $C_k^{[j]}$ to similarly define a map
$g_k^{[j]}\:X_k\to X_k$ with finitely many nontrivial point-inverses, all
of them shrinkable and of diameter $<\eps(k)$, and such that
$\acute p^j_k=g_k^{[j]}p^j_k$.

Let $\phi\:T_l\to\acute T_l$ be obtained by taking the quotient of
each $X_{l(i-1)}\x [i-\frac32,i-1]\i T_l$, identified with
$\Cyl(g_{l(i-1)}^{[l(i)]})$, by the projection
$\Cyl(g_{l(i-1)}^{[l(i)]})\to X_{l(i-1)}$ onto the image of
$g_{l(i-1)}^{[l(i)]}$.
Since $C_{l(i-1)}^{[l(i)]}$ is disjoint from $A_i$,
$\phi'(S/T\cap\bar X_l^\del)=S/T\cap\bar X_l^\del\i\Sigma(\acute s^l_{1^+})$,
where $\phi'\:\bar T_l\but\Sigma(\phi)\to\bar{\acute T_l}$ lifts to
the restriction of $\phi\x\phi$.
We may assume that $S/T$ is disjoint from $\Sigma(\phi)$ and more generally
$\phi^{-1}(\phi(x))=\{x\}$ whenever $(x,y)\in S$.
Since the cone of $\phi$ is PL homeomorphic to a (genuine) cone, so is
that of $\phi'$, and consequently that of its restriction
$\psi\:\bar T_l\but(\Sigma(\phi)\cup S/T)\to\bar{\acute T_l}\but\phi'(S/T)$.
Hence $\psi$ is a homology equivalence, and so $\tta_\del(\acute T_l)$ can be
represented by a cocycle with support in $\acute S/T:=\phi'(S/T)$.
By the above, $\acute S/T\cap\bar X_l^\del$ is contained in
$\Sigma(\acute s^l_{1^+})$.

{\it In the remainder of the proof we shall not use the original $p^i_j$,
$T_l$, $s^l_h$, $S$; these symbols will be recycled to denote what was
previously referred to as $\acute p^i_j$, $\acute T_l$, $\acute s^l_h$,
$\acute S$.}

{\bf Step 3: $S/T\i\Sigma(s^l_1)$ after further amending each $p^{i+1}_i$ by
a self-map of $X_i$, $\eps(i)$-close to the identity.}
Let $S_i\i X_{l(i-1)}\x X_{l(i-1)}$ denote the projection of
$\Cl{S\cap\tl X_{l(i)}\x(i-1,i]}$, and consider
$R_i:=p^{l(i-1)}_i\x p^{l(i-1)}_i(S_i)$.
From Step 2, $p^{l(i-1)}_i\x p^{l(i-1)}_i$ sends both
$S\cap\tl X_{l(i-1)}$ and $p^{l(i)}_{l(i-1)}(S\cap\tl X_{l(i)})$ into
$\Delta_{X_i}$.
Since $n>2$, we may assume that $R_i$ projects injectively onto the first
factor of $X_i\x X_i$.
Let $R_i^\#$ be the preimage of $\alpha X_i^{(n-1)}$ under this projection,
which without loss of generality is a finite set.
Since $S\i\Sigma_\eps(s^l_1)$, for each $(x,y)\in R_i$ we have
$\dist(x,y)<\eps(i)$.
Let $P_i$ and $P_i^\#$ be the quotients of $\Cyl(R_i\to R_i/T)$ and
$\Cyl(R_i^\#\to R_i^\#/T)$ obtained by shrinking to points the point-inverses
of the projection of $(R_i\cap\Delta_{X_i})\x I$ onto the first factor;
we may assume that $R_i\cap\Delta_{X_i}\i R_i^\#$.
Since $\eps(i)<\mesh(\alpha X_i)$, the projection $R_i\i X_i\x X_i\to X_i$
extends to a map $\mu_i\:P_i\to X_i$ which sends $R_i^\#/T$ to
$\alpha^*X_i^{(2)}\cap\alpha X_i^{(n-2)}$ and $R_i/T$ to
$\alpha^*X_i^{(2)}\cap\alpha X_i^{(n-1)}$, and embeds $P_i^\#\but (R_i^\#/T)$
into $\alpha X_i^{(n-1)}\but\alpha X_i^{(n-2)}$ and
$P_i\but (R_i/T\cup P_i^\#)$ into $X_i\but\alpha X_i^{(n-1)}$.
Then $\mu_i(P_i^\#)$ is a disjoint union of shrinkable $1$-polyhedra
$D_{ij}$ of diameter $<\eps(i)$ each.
We may view $D_i:=\bigcup D_{ij}\x D_{ij}$ (or rather $D_i\cup\Delta_{X_i}$)
as an equivalence relation on $X_i$.
Now $\mu_i(P_i)/D_i$ is a union of $2$-polyhedra $E_{ij}$, disjoint from each
other except at the finite set $\alpha^*X_i^{(2)}\cap\alpha X_i^{(n-2)}$;
each $E_{ij}$ is of diameter $<\eps(i)$ and meets precisely one
open $(n-1)$-simplex of $\alpha X_i^{(n-1)}$.
The projection $P_i\to R_i/T$ descends to maps
$\nu_{ij}\:E_{ij}\to E_{ij}\cap\alpha X_i^{(n-1)}$, and it is easy to see
that if $N$ is the second derived neighborhood of $E_{ij}$ in $\alpha X_i$
relative to $\alpha X_i^{(n-2)}$, then $(N,\Fr N)$ is PL homeomorphic to
$(N/\nu_{ij},\Fr N)$ keeping $N\cap\alpha X_i^{(n-2)}$ fixed.%
\footnote{If $f\:A\supset B\to C$ is a partial map, under $A/f$ we mean
$A/(f\x f)^{-1}(\Delta_C)$.}
It follows that if $\nu_i\:\mu_i(P_i)\to\mu_i(R_i/T)$ lifts to the projection
$P_i\to R_i/T$, the quotient $\eps$-map $\Pi_i\:X_i\to X_i/\nu_i$,
shrinking the point-inverses of $\nu_i$ to points, is arbitrarily closely
approximable by a homeomorphism $H_i\:X_i\to X_i/\nu_i$, and
the self-map $f_i:=H_i^{-1}\Pi_i$ of $X_i$ is $\eps(i)$-close to
the identity.

Let us now write $\grave p^{i+1}_i$ for $f_i p^{i+1}_i$.
In the obvious fashion one defines $\grave p^i_j$, $\grave T_h$,
$\grave s^l_1\:\grave T_l\to\grave T_1$, $\eta\:T_l\to\grave T_l$
and $\eta'\:T_l\but\Sigma(\eta)\to\grave T_l$, and verifies that
$\tta_\del(\grave T_l)$ can be represented by a cocycle with support in
$\eta'(S)$, which is contained in $\Sigma(\grave s^l_1)$. \qed
\enddemo

\demo{(b)} Since the cocycle representing $\tta(T_l)$ is of codimension one,
its support can be chosen in a $1$-polyhedron $S/T\i\Sigma(s^l_1)$, where
$S\cup\Delta_{T_l}$ is a subpolyhedron of $T_l\x T_l$, and $S\i\tl T_l^\del$
is disjoint from the $(2n-1)$-skeleton
$\bigcup_i(\alpha X_{l(i)}\x\alpha X_{l(i)})^{(2n-2)}\x (i-1,i]$ for some
triangulations $\alpha X_i$ of $X_i$ such that each $p^{i+1}_i$ is simplicial
in $\alpha X_{i+1}$ and some subdivision of $\alpha X_i$.

The closure of $S$ in $T_l\x T_l$ may be viewed as the union of compact
$1$-polyhedra $S_i:=\Cl{S\cap(\tl X_{l(i)}\x(i-1,i]})$ in
the mapping cylinders of $r^i_{i-1}:=p^{l(i)}_{l(i-1)}\x p^{l(i)}_{l(i-1)}$.
Let $R_i^{[i-1]}$ be the projection of $S_i$ to
$X_{l(i-1)}\x X_{l(i-1)}$.
Set $R_i^{[i]}=(r^i_{i-1})^{-1}(P_i^{[i-1]})$ and
$R_i^{[j]}=r^{i-1}_j(P_i^{[1]})$ for $j<i-1$.
Let $i'$ be the maximal index such that $l(i')\le i$.
Then $R_i^{[i']}\i\Delta_{X_{l(i')}}$ since $S/T\i\Sigma(s^l_1)$.
Consider
$$R_i=R_i^{[i]}\x(i-1,i]\cup R_i^{[i-1]}\x (i-2,i-1]\cup\dots\cup
R_i^{[i'+1]}\x (i',i'+1].$$
Finally, the $2$-polyhedron $R:=\bigcup_i R_i\but\Delta_{T_l}$ is our
``shadow'' of $S$.
The transitive closure $Q$ of $R$, i.e\. the set of all pairs
$(x_1,x_r)\in\tl T_l^\del\i\tl T_l$ such that $(x_i,x_{i+1})\in R$ for
$i=1,\dots,r-1$, is a $2$-polyhedron contained in the transitive set
$\Sigma(s^k_1)$.

View $Q$ and $Q_i:=Q\cap\tl X_{l(i)}$ (or rather their unions with
$\Delta_{T_l}$ and $\Delta_{X_{l(i)}}$) as equivalence relations on $T_l$
and $X_{l(i)}$.
Then $s^l_1$ factors through $K:=T_l/Q$, which is the mapping
telescope of the inverse sequence formed by $K_i:=X_{l(i)}/Q_i$ and
the quotients $q^i_j$ of the bonding maps $p^{l(i)}_{l(j)}$.
This yields a surjection $q\:\bar T_l^\del\but(Q/T)\to\bar K^\del$, whose
$\Sigma(q)$ meets every level $\del^{-1}(t)$ in an $(n+1)$-polyhedron, and
therefore is itself an $(n+2)$-polyhedron.
Hence the cone of $q$ collapses onto an $(n+3)$-polyhedron, and therefore
$q$ induces an isomorphism on $2n$-cohomology, since $n\ge 4$.
The image of $\tta_\del(K)$ under the latter coincides by naturality of
the Euler class with the image of $\tta_\del(T_l)$ under the restriction
$H^{2n}(\bar T_l^\del)\to H^{2n}(\bar T_l^\del\but (Q/T))$, which is zero
since $S\i R\i Q$.
Thus $\tta_\del(K)=0$. \qed
\enddemo

\demo{(c)} For a sufficiently rapidly decreasing $\eps$, the inverse limit
of the modified $\{X_i,p^i_j\}$ is $X$ by a well-known lemma of M. Brown
\cite{Bro; Theorem 2}.
Now $\{K_i,q^i_j\}$ can be obtained from $\{X_i,p^i_j\}$ by a thinning out
of indices, the inverse operation, and a further thinning out of indices.
\qed
\enddemo

The following Corollary to Theorem 2.2 seems to be of more interest in
connection with the proof of this theorem, rather than with its applications.
Let us call a compactum $X$ {\it pseudo-embeddable} into $\R^m$ if
it is the limit of an inverse sequence of compact polyhedra $X_i$ and
PL maps between them such that each finite mapping telescope
$\Tel(X_k\to\dots\to X_1\to X_0)$ admits a level-preserving embedding
into $\R^m\x [0,k]$.
(Recall from Criterion 3.2 that level-preserving embeddability of
the entire infinite telescope into $\R^m\x [0,\infty)$ is equivalent
to embeddability of $X$ into $\R^m$.)

The ANR from Example 2.12 is an example of a compactum, pseudo-embeddable
but not embeddable into $\R^{2n}$.
More generally, pseudo-embeddability is equivalent to quasi-embeddability
for any inverse limit of compact polyhedra and cell-like maps between them.
Indeed, given a polyhedron $X_{i+1}\i\R^m$ and a cell-like map
$p^{i+1}_i\:X_{i+1}\to X_i$, the quotient map $\R^m\to\R^m/p^{i+1}_i$
shrinking to points the point-inverses of $p^{i+1}_i$ is the final map of
a pseudo-isotopy of $\R^m$ (by Siebenmann's theorem, extended to dimension
$4$ by Quinn), which yields a level-preserving embedding
$\Cyl(p^{i+1}_i)\emb\R^m\x [i,i+1]$, extending the given embedding
$X_{i+1}\emb\R^m$.

In particular, every AR from Problem 1.5 pseudo-embeds into $\R^{2n}$.

\proclaim{Corollary 3.5} An $n$-dimensional compactum $X$, $n>3$, embeds
into $\R^{2n}$ if it pseudo-embeds into $\R^{2n-1}$.
\endproclaim

It follows that Akhmetiev's compactum from Example 2.15, which quasi-embeds
into $\R^{2n-1}$, does not pseudo-embed into $\R^{2n-1}$.

\demo{Proof} Let $T_\infty$ be the infinite mapping telescope of the given
inverse sequence.
By the hypothesis, every finite initial telescope $T_i\i T_\infty$ embeds
into $\R^{2n-1}$.
Hence $w_1(\tl T_i^\del\to\bar T_i^\del)^{2n-1}=0$ for each $i$.
Since the groups $H^{2n-2}(\bar T_i^\del;\Z/2)$ are finite, their inverse
sequence satisfies the Mittag-Leffler condition.
Then from the Milnor sequence,
$w_1(\tl T_\infty^\del\to\bar T_\infty^\del)^{2n-1}=0$.
Hence its Bockstein image (cf\. \cite{M5})
$\tta^\del(T_\infty)=e(\tl T_\infty^\del\to\bar T_\infty^\del)^{2n}$ vanishes
as well.
By Lemma 3.3 and Criterion 3.2, $X$ embeds into $\R^{2n}$.
(Alternatively, from the proof of Lemma 2.16 one
can deduce that $\tta(X)=0$ and refer directly to Theorem 2.2.) \qed
\enddemo

The following lemma was used to deduce Corollary 2.6.

\proclaim{Lemma 3.6} (a) If $H^n(X,X\but x)\to H^n(X)$ is onto for
some $x\in X$, $\hat H^{2n}(\bar X)=0$.

(b) If in addition $H^{n+1}(X,X\but x)=0$ for every $x\in X$,
$H^{2n}(\bar X)=0$.
\endproclaim

In particular, $H^n(X\but x)=0$ for every $x\in X$ implies
$H^{2n}(\bar X)=0$.
Note also that $H^{n+1}(X,X\but x)$ is either zero or uncountable
\cite{Harl}.

\demo{Proof} The hypothesis of (a) implies that $H^n(X)$ lies in the image of
$\invlim H^n_c(U_i)$, where $U_i$ is a fundamental sequence of open
neighborhoods of the given point $x$.
By the functoriality in the K\"unneth formula \cite{Bre},
$H^{2n}(X^2)$ lies in the image of $\invlim H^{2n}_c(U_i^2)$.
So $H^{2n}(X^2,X^2\but (x,x))\to H^{2n}(X^2)$ is onto.
Hence $H^{2n}(K)=0$ for every compact $K\i\tl X$.
Then by one of the Smith sequences, $H^{2n}(C)=0$ for every compact
$C\i\bar X$, which completes the proof of (a).
To deduce that $H^{2n}(\tl X)=0$ (which will similarly imply (b)), it remains
to show that $H^{2n+1}(X^2,\tl X)=0$.

Let $x\in X$ and let $U_i(x)$ be a fundamental nested sequence of open
neighborhoods of $x$.
By the Milnor exact sequence,
$\lim^1 H^n(X,X\but U_i(x))=H^{n+1}(X,X\but x)$, which is zero by
the hypothesis.
Since each $H^n(X,X\but U_i(x))$ is countable (being a direct limit of
cohomology of compact polyhedral pairs), the vanishing of the derived limit
is equivalent to the Mittag-Leffler condition \cite{Gr}.
For convenience of notation we shall use cohomology with compact support
$H^n_c(U_i(x))=H^n(X,X\but U_i(x))$.
Thus for each $p$ and each $x\in X$ there exists a $q_x\ge p$ such that for
each $r\ge q_x$ the image of $H^n_c(U_r(x))$ in $H^n_c(U_p(x))$ does not
depend on $r$.
By the functoriality in the K\"unneth formula \cite{Bre}, so does the image
of $H^{2n}_c(U_r(x)^2)$ in $H^{2n}_c(U_p(x)^2)$.

Consider the fundamental nested sequence
$V_i:=\bigcup_{x\in X}U_i(x)^2$ of open neighborhoods of $\Delta_X$.
(The use of the uncountable Axiom of Choice can be eliminated using
compactness of $X$.)
Since $X$ is compact, there exists a $q\ge p$ such that $V_q$ is contained
in $V:=\bigcup_{x\in S} U_{q_x}(x)^2$ for some finite $S\i X$;
without loss of generality, $q\ge q_x$ for each $x\in S$.
By an iterated application of the Mayer--Vietoris sequence,
$\bigoplus_{x\in S} H^{2n}_c(U_{q_x}(x)^2)\to H^{2n}_c(V)$
is onto.
Therefore if $r\ge q$,
$H^{2n}_c(\bigcup_{x\in S} U_r(x)^2)\to H^{2n}_c(V)$ is onto.
Since $\bigcup_{x\in S} U_r(x)^2\i V_r\i V\i V_p$, the image of
$H^{2n}_c(V_r)$ in $H^{2n}_c(V_p)$ does not depend on $r$, as long as
$r\ge q$.
Thus $H^{2n+1}(X^2,\tl X)=\lim^1 H^{2n}_c(V_i)=0$. \qed
\enddemo

\head 4. Results in the metastable range \endhead

The definition of the extraordinary van Kampen obstruction
$$\Theta^m(X)\in\omega^{mT}_{\Z/2}(\tl X_+)\simeq [\tl X_+,S^{mT}]_{\Z/2}$$
to embeddability of $X$ into $\R^m$ as given in \cite{M5} makes perfect sense
if $X$ is a compactum of dimension $<m$.
Here $K_+$ denotes the pointed $\Z/2$-space $K\sqcup *$, and
$S^{mT}$ is the $m$-sphere with the action of $\Z/2$ fixing the basepoint $*$
and restricting to the sign action $x\inv -x$ on the complement of $*$,
identified with $\R^m$.
The obstruction $\Theta^m(X)$ is the class of the composition
$$\tl X@>>> S^\infty@>>> S^\infty/S^{m-1}=S^\infty_+\wedge S^{mT}@>>>S^{mT}$$
of an arbitrary equivariant map $\tl X\to S^\infty$ (with respect to the
antipodal action of $\Z/2$ on $S^\infty$) and the obvious projections.
If $X$ embeds in $\R^m$, such a map factors up to equivariant homotopy
through $S^{m-1}\i S^\infty$, hence $\Theta^m(X)=0$.
Let $\hat\Theta^m(X)$ be the image of $\Theta^m(X)$ in
$\hat\omega^{mT}_{\Z/2}(\tl X_+):=\lim\omega^{mT}_{\Z/2}(K_+)$, the inverse
limit over all compact invariant subsets $K\i\tl X$.

\proclaim{Lemma 4.1} \cite{M5; proof of Lemma 4.4} $\Theta^m(X)=0$ iff there
exists an equivariant map $\tl X\to S^{m-1}$, and $\hat\Theta^m(X)=0$ iff
such a map exists over each compact invariant subset $K\i\tl X$.
\endproclaim

Lemma 4.1 will {\it not} be used in the proof of

\proclaim{Theorem 4.2} Let $X$ be an $n$-dimensional compactum.

(a) $X$ embeds into $\R^m$, $m>\frac{3(n+1)}2$, if and only if
$\Theta^m(X)=0\in\omega^{mT}_{\Z/2}(\tl X_+)$.

(b) $X$ quasi-embeds into $\R^m$, $m\ge\frac{3(n+1)}2$, if and only if
$\hat\Theta^m(X)=0\in\hat\omega^{mT}_{\Z/2}(\tl X_+)$.
\endproclaim

\demo{Proof}
Modulo its polyhedral case, which is proved in \cite{M5}, Theorem 4.2 follows
by the same argument as Theorems 2.2 and 2.4, using that $\omega^*_{\Z/2}$ is
a generalized $\Z/2$-equivariant cohomology theory, cf\. \cite{M5}.
The grading is ``upside down'' in that $\omega^{mT}_{\Z/2}(Y)$ vanishes for
$m>\dim Y$, but may be nontrivial for negative $m$.
In particular, the vanishing of $H^{2n}(P^k/(\Z/2))$, where $2n>k$, now
becomes that of $\omega^{mT}_{\Z/2}(P^k)$, where $m>k$.
(This was used in the very ends of the proofs of Propositions 3.1 and 3.4b,
where the cones of certain maps, denoted by $q$ in both cases, collapse
onto a $k$-polyhedron $P$, where $k=n+(2n-m)+1$ and $n+(2n-m+1)+1+1$,
respectively.)
With these matters taken into account, only straightforward changes are to
be made in the proofs, including a routine generalization of steps 2 and 3
(which involved no cohomology whatsoever) in the proof of Proposition 3.4a,
and a simplification of the proof of Lemma 3.3, where the computation of
the second obstruction is no longer necessary. \qed
\enddemo

\proclaim{Proposition 4.3} Let $X$ be an $n$-dimensional compactum and $Y$
a compactum of dimension $\le\frac k2+\frac{m-n-3}2$ such that there exists
an equivariant map $S^{k-1}\to\widetilde{CY}$.
Then $\Theta^{m+k}(X*Y)=0$ implies $\Theta^m(X)=0$.
\endproclaim

\demo{Proof} The intersection of $\widetilde{CY}$ with $Y\x CY\cup CY\x Y$
is equivariantly homotopy equivalent to the double mapping cylinder $Z$ of
the projections $Y@<p_1<<\tl Y@>p_2>>Y$.
An equivariant retraction of $\widetilde{CY}$ onto this intersection is
given by $((x,t),(y,s))\mapsto((x,\max\{t-s,0\}),(y,\max\{s-t,0\}))$, where
$CY=Y\x I/Y\x\{1\}$.
(In fact, this is the final map of an equivariant deformation retraction.)
Thus there exists an equivariant map $S^{k-1}\to Z$.

$\widetilde{X*Y}$ is equivariantly homotopy equivalent to its subset $S$
that is the quotient of $\tl X\x\tl Y\x I\x I$.
Shrinking to points $X\x\{(y,y',1,t)\}$ and $X\x\{(y,y',t,1)\}$ for each
$(y,y')\in\tl Y$ and $0<t<1$ as well as $X\x\{(y,1,0)\}$ and $X\x\{(y,0,1)\}$
for each $y\in Y$ yields an equivariant map $S\to\tl X*Z$.
The relative mapping cylinder of this map equivariantly collapses onto a pair
that is excision-equivalent to $(CX,X)\x Z/\phi$, where the partial map
$\phi\:(CX,X)\x Z\supset (CX,X)\x\tl Y\to (I,\{0\})\x\tl Y$ shrinks each
level of the cone.
Since this pair is at most $(m+k-1)$-dimensional,
$\omega^{(m+k)T}_{\Z/2}(\widetilde{X*Y}_+)\simeq
\omega^{(m+k)T}_{\Z/2}(\tl X*Z_+)$.
Therefore the original equivariant map $\widetilde{X*Y}\to S^{m+k-1}$ yields
an equivariant map $\Sigma^k\tl X\to Z*\tl X\to S^{m+k-1}$.
Thus $\Theta^m(X)=0$. \qed
\enddemo

\proclaim{Corollary 4.4} Let $X$ be an $n$-dimensional compactum,
$m>\frac{3(n+1)}2$, and $k>0$.
The following assertions are equivalent:

(i) $X$ embeds into $\R^m$;

(ii) the $k$-fold suspension $\Sigma^k X$ embeds into $\R^{m+k}$;

(iii) the $k$-fold cone $C^k X$ embeds into $\R^{m+k}$;

(iv) $X*T_k$ embeds into $\R^{m+2k}$, where $T_k$ is the join of $k$ copies
of $\{1,2,3\}$;

(v) $X*Z_k$ embeds into $\R^{m+2k}$, where $Z_k$ is the $(k-1)$-skeleton of
the $2k$-simplex.
\endproclaim

\demo{Proof} (ii)$\imp$(i) and (iii)$\imp$(i) follow by a {\it repeated}
application of Proposition 4.3.
(iv)$\imp$(i) and (v)$\imp$(i) follow since $\widetilde{CY_k}$ and
$\widetilde{CZ_k}$ contain $S^{2k-1}$ by the Flores construction
(see \cite{M5}). \qed
\enddemo

A shorter proof of (ii)$\imp$(i) is given by

\proclaim{Theorem 4.5} Let $X$ be an $n$-dimensional compactum,
$Y$ a compactum such that there exists an equivariant map $S^{k-1}\to\tl Y$,
and assume $m>\frac{3(n+1)}2$.
If $X*Y$ admits a level-preserving embedding into $\R^{m+k}\x I$, then $X$
embeds into $\R^m$.
\endproclaim

\demo{Proof} The subset of $\widetilde{X*Y}$ consisting of all pairs $(a,b)$
such that $p(a)=p(b)$, where $p\:X*Y\to I$ is the projection, is homeomorphic
to $\tl X*\tl Y$.
So the given embedding yields an equivariant map
$\Sigma^k\tl X\to\tl X*\tl Y\to S^{m+k-1}$, thus $\Theta^m(X)=0$. \qed
\enddemo

\proclaim{Proposition 4.6} Let $X$ be an acyclic $n$-dimensional compactum
and $Y$ a compactum of dimension at most $\frac k2+\frac{m-2}2$ such that
there exists an equivariant map $S^{k-1}\to\widetilde{CY}$.
If $\Theta^{m+k}(X\x CY)=0$, then $\Theta^m(X)=0$.
\endproclaim

\demo{Proof} Let $S$ be the subset of $\widetilde{X\x CY}$ that is
the quotient of
$\tl X\x\tl Y\x I\x I\cup X\x X\x\tl Y\x (I\x\{1\}\cup\{1\}\x I)$.
Shrinking to points $X\x X\x\{(y,y',t,1)\}$ and $X\x X\x\{(y,y',1,t)\}$
for each $(y,y')\in\tl Y$ as well as $X\x X\x\{(y,0,1)\}$ and
$X\x X\x\{(y,1,0)\}$ for each $y\in\tl Y$ yields an equivariant map
$S\to\tl X*Z$, where $Z$ is the double mapping cylinder of the projections
$Y@<p_1<<\tl Y@>p_2>>Y$.
The relative mapping cylinder of this map equivariantly collapses onto a pair
that is excision-equivalent to $Z\x(C(X\x X),X\x X)$.
Since $X$ is acyclic, $\Sigma(X\x X)$ is contractible, so the projection
$Z\x\Sigma(X\x X)\to Z$ is an equivariant homotopy equivalence.
Hence
$\omega^{iT}_{\Z/2}(Z\x(C(X\x X),X\x X)_+)\simeq\omega^{iT}_{\Z/2}(Z_+)$.
Therefore, since $Z$ is at most $(m+k-1)$-dimensional, we get
$\omega^{(m+k)T}_{\Z/2}(S_+)\simeq\omega^{(m+k)T}_{\Z/2}(\tl X*Z_+)$.
Since there exists an equivariant map $S^{k-1}\to Z$ (see the proof of
Proposition 4.3), the given equivariant map
$S\i\widetilde{X\x CY}\to S^{m+k-1}$ gives rise to an equivariant map
$\Sigma^k\tl X\to Z*\tl X\to S^{m+k-1}$.
Thus $\Theta^m(X)=0$. \qed
\enddemo

\proclaim{Corollary 4.7} Let $X$ be an acyclic $n$-dimensional compactum,
$m>\frac{3(n+1)}2$, $k>0$.
The following assertions are equivalent:

(i) $X$ embeds into $\R^m$;

(ii) $X\x I^k$ embeds into $\R^{m+k}$;

(iii) $X\x \text{(triod)}^k$ embeds into $\R^{m+2k}$.
\endproclaim

The hypothesis of acyclicity can be weakened.
The proof of (ii)$\imp$(i) works for all $X$ such that the (already) stable
cohomotopy groups $\omega^{m+i}(X\x X)=0$ for $1\le i\le k$.
The proof of (iii)$\imp$(i) works for all $X$ such that
$\omega^{(m+2i)T}_{\Z/2}(X\x X\x S^1_+)=0$ for $1\le i\le k$, where $\Z/2$
acts by exchanging the first two factors and transforming $S^1$ by
the antipodal involution.

On the other hand, the dimensional restrictions cannot be dropped.
If $X$ is a non-simply-connected homology $n$-ball (i.e\. a homology sphere
minus an open ball), then $X$ does not embed in $\R^n$, but $X\x I$ embeds
in $\R^{n+1}$ since every homology sphere bounds a contractible topological
manifold (Kervaire, Freedman--Quinn), whose double has to be the sphere by
Seifert--van Kampen and the generalized Poincar\'e conjecture.
One can increase the codimension in this example by considering appropriate
spines of $X$.

A shorter proof of (ii)$\imp$(i) is given by

\proclaim{Theorem 4.8} Let $X$ be an acyclic $n$-dimensional compactum,
$m>\frac{3(n+1)}2$.
If $X$ can be instantaneously taken off itself in $\R^{m+1}$, i.e\.
the mapping cylinder of the projection $X\sqcup X\to X$ admits
a level-preserving embedding into $\R^{m+1}\x I$, then $X$ embeds in $\R^m$.
\endproclaim

In fact, $\theta(X)$ may be thought of as precisely the obstruction to
isotopic realizability (see \cite{M1}) of the composition of
the projection $X\sqcup X\to X$ and some embedding $X\emb\R^{2n+1}$.

\demo{Proof}
Let $g\:X\emb\R^{m+1}$ be the given embedding that can be instantaneously
taken off itself.
Then $\tl g\:\tl X\to S^m$ is homotopic to a map that extends to $X\x X$.
Since $X$ is acyclic, this map must be null-homotopic.
Combining the null-homotopy with its reflection, we get an equivariant map
$\Sigma\tl X\to S^m$.
Hence $\Theta^m(X)=0$. \qed
\enddemo

\remark{Remark} F. Quinn has shown that every tame embedding of an ANR into
a manifold of dimension $>4$ has a mapping cylinder neighborhood
(see \cite{Q}).
By an argument of M. Cohen, this implies that $X\x I^{2n+1}$ is collapsible
for every $n$-dimensional AR $X$, $n>1$.
Indeed, by \cite{Sht}, $X$ tamely embeds in $\R^{2n+1}$, where its Quinn's
neighborhood is homeomorphic to $I^{2n+1}$.
For it is contractible and has a simply-connected boundary, since $X$ is
tame and so misses generic $2$-disks.
Let $f\:S^{2n}\to X$ be the map such that $(I^{2n+1},X)=(\Cyl(f),X)$.
Now $X\x I^{2n+1}$ is the mapping cylinder of the projection
$\pi\:X\x S^{2n}\to X$.
It collapses onto $\Cyl(\pi|_{\Gamma_f})$, where $\Gamma_f\i S^{2n}\x X$
is the graph of $f$.
But this is homeomorphic to $I^{2n+1}$, which is collapsible.
\endremark

\proclaim{Corollary 4.9} An acyclic $n$-dimensional compactum, $n>3$, embeds
in $\R^{2n}$ if it immerses there.
\endproclaim

\demo{Proof} Let $\phi\:X\imm\R^{2n}$ be an immersion.
We may think of it as a composition
$X\overset g\to\emb M\overset\psi\to\imm\R^{2n}$, where $M$ is a
PL $2n$-manifold.
Since $n>2$, we may assume that $g$ is tame \cite{Sht}.
This means that for each $\eps>0$, there exists an $\eps$-pseudo-isotopy,
taking $g(X)$ onto a subpolyhedron of $M$.
It follows that $g$ can be replaced by an embedding $g'$ such that
the immersion $\phi'=\psi g'$ has at most $0$-dimensional double point set
$\Delta_{\phi'}=(\phi'\x\phi')^{-1}(\Delta_{\R^{2n}})\but\Delta_X\i\tl X$
and no triple points (compare the arguments in \cite{SS}).

Since $\Delta_{\phi'}$ is $0$-dimensional, there exists an equivariant map
$\Delta_{\phi'}\to S^0$.
Extending it to a continuous function $X\to\R$, where $\Delta_{\phi'}$ is
identified with its homeomorphic projection to the first factor of $X\x X$,
and $S^0\i\R$, we obtain that the immersion $\phi'$ lifts vertically to
an embedding $\hat\phi'\:X\to\R^{2n+1}$.
Then the translates of $\hat\phi'(X)$ by sufficiently small distances in
the vertical direction are all disjoint from each other, which yields
an embedding $X\x I\emb\R^{2n+1}$.
Since $n>3$, by Corollary 4.7 $X$ embeds in $\R^{2n}$. \qed
\enddemo

Taras Banakh asked one of the authors in July 2006: what is the minimal
dimension of Euclidean space containing $\mu^n\x I^k$, where $\mu^n$ is
the universal Menger compactum?
Since there exist contractible $n$-dimensional compacta, non-embeddable into
$\R^{2n}$ (see Examples 1.2, 2.17), we obtain from Corollary 4.7
the following answer%
\footnote{Non-embeddability of an $n$-dimensional compactum cross $I^k$ into
$\R^{2n+k}$ can also be proved directly.
The proof of Proposition 4.6 yields an equivariant homotopy equivalence
$\Sigma\tl P\simeq\widetilde{P\x I}$, and hence
$\Sigma^k\tl P\simeq\widetilde{P\x I^k}$ for every contractible
polyhedron $P$.
One only needs to plug it for $P=$ the product of $n$ triods (or $P=$ the
cone over the $(n-1)$-skeleton of the $2n$-simplex) into one of the proofs
(in \cite{RS}, \cite{KR}, \cite{RSS} or in Example 2.17 above) that
$P\x\N^*$ does not embed into $\R^{2n}$, to get that $P\x I^k\x \N^*$ does
not embed into $\R^{2n+k}$.}
to Banakh's question:

\proclaim{Corollary 4.10} $\mu^n\x I^k$ does not embed into $\R^{2n+k}$.
\endproclaim

The motivation for Banakh's question was as follows.
By general position $\mu^n\x I^k$ embeds in $\R^{2n+k+1}$.
If this estimate were not sharp for $k=2n$, one would deduce the following
result, answering \cite{BCZ; Question 1.4} (see also \cite{BV; Problem 6}).

\proclaim{Proposition 4.11} $\mu^{2n}$ does not embed into $\mu^n\x I^{2n}$.
\endproclaim

This result, answering the question from ``Open Problems in Topology II'', is
proved by the following self-contained argument.

\demo{Proof}
Let $P$ be the $2n$-skeleton of the $(4n+2)$-simplex, or any other
$2n$-polyhedron, non-embeddable in $\R^{4n}$.
If $\mu^{2n}$ embeds in $\mu^n\x I^{2n}$, so does $P$.
Hence for each $\eps>0$ it admits an $\eps$-map to $Q_\eps\x I^{2n}$ for some
$n$-polyhedron $Q_\eps$ (an appropriate nerve of the Menger cube).
Now $Q_\eps\x I$ embeds in $\R^{2n+1}$ \cite{RSS} since $Q_\eps$ immerses in
$\R^{2n}$ with isolated double points.
Hence $Q_\eps\x I^{2n}$ embeds in $\R^{4n}$, so $P$ quasi-embeds there.
But for polyhedra quasi-embeddability is equivalent to embeddability in
the case of double dimension (except $2\to 4$) by completeness of the van
Kampen obstruction (also it is in the metastable range by the
Haefliger--Weber criterion).
Thus $P$ embeds in $\R^{4n}$, contradicting its choice. \qed
\enddemo

We recall (see \cite{M5}) that the Conner--Floyd cohomological co-index of
$\tl X$ is the maximal $m$ such that the Hurewicz image
$h(\Theta^m(X))\in H^m(\bar X;\Z\T^{\otimes m})$ is nonzero, where
$\Z\T^{\otimes m}$ is induced from the orientation sheaf of $\RP^{m-1}$
under the classifying map $\bar X\to\RP^\infty$.

\proclaim{Theorem 4.12} Suppose that $m>\frac{3(n+1)}2$.

(a) All $n$-dimensional compacta $X$ with $H^i(\tl X)=0$ for $i\ge m$
embed in $\R^m$.

(b) Let $X$ be an $n$-dimensional locally acyclic compactum with
$\hat H^i(\tl X)=0$ for $i>m$.
Then $X$ embeds in $\R^m$ if and only if $\coind_\Z(\tl X)<m$.
\endproclaim

We note a certain similarity of (b) with the homological criterion for
isotopic realizability under the assumption of discrete realizability by
skeleta \cite{AM1}, \cite{M2}.

\demo{Proof} By Theorem 4.2 and the preceding remark, it suffices to show
existence of an equivariant map $\tl X\to S^{m-1}$.
By obstruction theory, such a map exists if the first obstruction
$h(\Theta^m(X))\in H^m(\bar X;\Z\T^{\otimes m})$ vanishes (which is
already given in (b)) and the higher obstruction groups
$H^{i+1}(\bar X;\pi_i(S^{m-1})\otimes\Z\T^{\otimes m}))$ are trivial.

Let us represent $\tl X$ as a union of compact subsets $K_j$ with each
$K_j\i\Int K_{j+1}$.
Fix some $i\ge m-1$ and $k\ge m-1$ for (a), $i\ge m$ and $k\ge m$ for (b).
By the Milnor exact sequence, the hypothesis of (a) implies
$\hat H^{i+1}(\tl X)=0$ and $\lim^1 H^i(K_j)=0$.
By the universal coefficients formula for compacta \cite{Bre},
$\hat H^{i+1}(\tl X;\pi_k(S^{m-1}))=0$ in both (a) and (b), and
$\lim^1 H^i(K_j;\pi_k(S^{m-1}))=0$ in (a).
On the other hand, since $X$ is locally acyclic in (b), and
$K_j\i\Int K_{j+1}$, the image of $H^i(K_{j+1};\pi_k(S^{m-1}))$ in
$H^i(K_j;\pi_k(S^{m-1}))$ is finitely generated for each $j$.
By Serre's theorem $\pi_k(S^{m-1})$ is finite, hence so is this image.
Thus the inverse sequence $H^i(K_j;\pi_k(S^{m-1}))$ satisfies the
Mittag-Leffler condition, so $\lim^1 H^i(K_j;\pi_k(S^{m-1}))=0$ in (b)
as well.

Now the Milnor sequence implies that $H^{i+1}(\tl X;\pi_k(S^{m-1}))=0$
for all $i,k\ge m-1$ in (a) and for all $i,k\ge m$ in (b).
It follows from the two Smith sequences, by downward induction on $i$, that
$H^{i+1}(\bar X;\pi_k(S^{m-1}))=0=H^{i+1}(\bar X;\pi_k(S^{m-1})\otimes\Z\T)$
for all $i,k\ge m-1$ in (a) and for all $i,k\ge m$ in (b).
In particular, this holds for $i=k$. \qed
\enddemo

\proclaim{Corollary 4.13} Every $n$-dimensional compactum $X$ with
$H^{n-d}(X\but x)=0$ for each $x\in X$ and $d\le k$, where $k<\frac{n-3}2$,
embeds into $\R^{2n-k}$.
\endproclaim

Another proof of the case $k=0$ is contained in the proof of Corollary 2.6
(see remark to Lemma 3.6).

\demo{Proof} Let us represent $\tl X$ as a union of compact subsets
$K_1\i K_2\i\dots$.
Let $K_i^x$ stand for the intersection of $K_i$ with $(X\but x)\x x$.
From Milnor's exact sequence, $\lim H^{n-d}(K_i^x)=0=\lim^1 H^{n-d-1}(K_i^x)$
for $d\le k$ (and each $x$).
Since each $H^{n-d}(K_i^x)$ is countable (being a direct limit of cohomology
of compact polyhedra), each inverse sequence $H^{n-d}(K_i^x)$ satisfies
the Mittag-Leffler condition \cite{Gr}.
On the other hand, when $d\le k$, this sequence has trivial inverse limit,
that is, each $c\in H^{n-d}(K_i^x)$ has a nonempty preimage in
$H^{n-d}(K_j^x)$ for only finitely many $j>i$.
Combined with the Mittag-Leffler condition, this implies that for each $i$
there exists a $j>i$ such that the image of $H^{n-d}(K_j^x)$ in
$H^{n-d}(K_i^x)$ is zero.

Let $\H^{n-d}(\pi_i)$ be the Leray sheaf of the projection
$\pi_i\:K_i\i X\x X\to X$ \cite{Bre}.
Since $K_i$ is compact, $\pi_i$ is closed, hence the stalks
$\H^{n-d}(\pi_i)_x\simeq H^{n-d}(K_i^x)$.
\footnote{Note that the projection $\pi\:\tl X\i X\x X\to X$ is not closed,
and that $\H^{n-d}(\pi)_x$, which is by definition the direct limit of
$H^{n-d}(\tl X\cap X\x U_i(x))$, where $U_i(x)$ is a fundamental nested
sequence of neighborhoods of $x$, may in general differ from
$H^{n-d}(X\but x)$. This is why we need the $K_i$'s.}
Consider the Leray spectral sequence \cite{Bre}
$$H^p(X;\H^q(\pi_i))\quad\Rightarrow\quad H^{p+q}(K_i).$$
Since $\H^{n-d}(\pi_j)\to\H^{n-d}(\pi_i)$ is zero for an appropriate $j>i$
and all $d\le k$, so is $H^{2n-d}(K_j)\to H^{2n-d}(K_i)$.
Also, since $\H^{n-k-1}(\pi_j)\to\H^{n-k-1}(\pi_i)$ does not depend on $j>j_0$
for an appropriate $j_0>i$, so does $H^{2n-k-1}(K_j)\to H^{2n-k-1}(K_i)$.
Hence the sequence $H^{2n-d}(K_i)$ has zero inverse limit when $d\le k$
and zero derived limit when $d\le k+1$.
It now follows from the Milnor sequence that $H^{2n-d}(\tl X)=0$ for all
$d\le k$. \qed
\enddemo

Note that the hypothesis of Corollary 4.13 is satisfied if

\noindent
(i) $H^{n-d}(X)\simeq H^{n-d}(S^n)$ for $0\le d\le k$,

\noindent
(ii) $H^{n-d}(X,X\but x)\simeq H^{n-d}(\R^n,\R^n\but\{0\})$
for $-1\le d\le k-1$ and each $x\in X$, and

\noindent
(iii) $H^n(X,X\but x)\to H^n(X)$ is an isomorphism for each $x\in X$.

So, since generalized manifolds satisfy Poincar\'e duality, we
obtain an elementary proof of the following result
(using Theorem 1.4 in the case $k=0$, where no orientability is assumed).

\proclaim{Corollary 4.14 (Bryant--Mio)} \cite{BM}
Every homologically $k$-connected $n$-dim\-ensional generalized manifold,
$k<\frac{n-3}2$, embeds in $\R^{2n-k}$.
\endproclaim

V. M. Buchstaber asked the first author in September 2005, whether there is
a generalization to polyhedra of the classical Penrose--Whitehead--Zeeman
Theorem that $k$-connected manifolds embed in $\R^{2n-k}$ in the metastable
range. (See \cite{Z}, which includes Irwin's extension to codimension three.)
Corollary 4.13 gives an answer, which however is not as revealing as one
might expect.

Indeed, $n$-polyhedra satisfying condition (ii) have links of $p$-simplices,
$p\ge 0$, with top $k-1$ cohomology groups isomorphic to those of
$S^{n-p-1}$.
When $k\ge 1$, such polyhedra are polyhedral homology manifolds with
codimension $k+2$ singularities.
For these, Poincar\'e duality holds in the first $k$ dimensions, so condition
(i) along with (iii) is now equivalent to homological $k$-connectedness:
$\tl H^d(X)=0$ for $d\le k$.
Thus $n$-polyhedra, satisfying (i)--(iii) with $k\ge 1$ are homologically
$k$-connected polyhedral homology manifolds with (restricted) codimension
$k+2$ singularities.

At the same time, the original PWZ argument works for $k$-connected (genuine)
manifolds with arbitrary codimension $k+2$ singularities, as all
constructions in this argument are disjoint from such singularities by
general position.
As for the case $k=0$, Sarkaria \cite{Sa} noticed that the PWZ method works
to embed in $\R^{2n}$ every quotient $X$ of a PL $n$-manifold $M$, $n>2$, by
a PL identification on the boundary such that no two components of $M$ remain
disjoint in $X$.
Clearly, this includes all $n$-polyhedra satisfying the hypothesis of
Theorem 1.3.

\enddocument
\end